# IMAGINARY QUADRATIC FIELDS WITH $\mathrm{Cl}_2(k) \simeq (2,2,2)$

E. BENJAMIN
F. LEMMERMEYER
C. SNYDER

ABSTRACT. We characterize those imaginary quadratic number fields, $k$, with 2-class group of type $(2,2,2)$ and with the 2-rank of the class group of its Hilbert 2-class field equal to 2. We then compute the length of the 2-class field tower of $k$.

## 1. Introduction

Let $K$ be an algebraic number field and $\mathrm{Cl}_p(K)$ the Sylow $p$-subgroup of its ideal class group, $\mathrm{Cl}(K)$. Denote the orders of $\mathrm{Cl}(K)$ and $\mathrm{Cl}_p(K)$ by $h(K)$ and $h_p(K)$, respectively. Let $K^1$ denote the Hilbert $p$-class field of $K$ (always in the wide sense). Finally, for nonnegative integers $n$, let $K^n$ be defined inductively as $K^0 = K$, $K^{n+1} = (K^n)^1$, and $K^\infty = \bigcup K^n$. Then the sequence

$$K^0 \subseteq K^1 \subseteq \cdots \subseteq K^n \subseteq \cdots \subseteq K^\infty$$

is called the $p$-class field tower of $K$. In general, very little is known about this tower; for instance, even its length has not been determined for most fields.

Now let $k$ be an imaginary quadratic number field and $p = 2$. Then by the work of Golod and Shafarevich [12], (also see [26]), the 2-class field tower of $k$ is infinite, if $\operatorname{rank} \mathrm{Cl}_2(k) \geq 5$. (Here rank means minimal number of generators.) When $\operatorname{rank} \mathrm{Cl}_2(k) = 2$ or $3$, there are examples of $k$ with infinite 2-class field towers as well as those for which the tower is finite. For $\operatorname{rank} \mathrm{Cl}_2(k) = 4$ no example of such a $k$ with finite 2-class field tower has ever been exhibited. In fact, it has been conjectured [24],[25] that the 2-class field tower of an imaginary quadratic field with $\operatorname{rank} \mathrm{Cl}_2(k) = 4$ is always infinite. On the other hand, until recently, all known examples of number fields with finite 2-class field tower (even more generally for $p$-class field towers) have length $\leq 2$. In [9] M. Bush exhibited an example of an imaginary quadratic field, namely $k = \mathbb{Q}(\sqrt{-445})$, for which the 2-class field tower is of length 3.

We consider another approach to finding a field $K$ with finite $p$-class field tower of length $\geq 3$, and at the same time continue our quest for a complete classification of complex quadratic number fields with 2-class field towers of length 2. We consider those $K$ with $\operatorname{rank} \mathrm{Cl}_p(K^1) = 2$. If $\operatorname{rank} \mathrm{Cl}_p(K) = 2$, then by a group theoretic result of Blackburn [7], Theorem 4, $\mathrm{Cl}_p(K^2)$ is trivial, whence the $p$-class field tower of $K$ has length 2. If, however, $\operatorname{rank} \mathrm{Cl}_p(K) \geq 3$, then by [7], Theorem 1,

The last author would like to thank the Department of Mathematics of the University of Maine for providing travel funds to the University of Cologne to work on the results in this and other papers during his sabbatical leave.

$\mathrm{Cl}_p(K^2)$ is either cyclic or trivial. In the nontrivial cyclic case, $K$ would have $p$-class field tower of length 3.

In this article we study the problem just mentioned when $k$ is an imaginary quadratic field and $p = 2$. We must have $\operatorname{rank}\mathrm{Cl}_2(k) \geq 3$; but if $\operatorname{rank}\mathrm{Cl}_2(k) \geq 4$, then it can be shown that $\operatorname{rank}\mathrm{Cl}_2(k^1) \geq 5$ (see for example the method of Propositon 3 of [2]). Hence $\operatorname{rank}\mathrm{Cl}_2(k) = 3$ is our only viable option if we want to show that the 2-class field tower has length 2. We assume that the rank is thus 3 and furthermore that $\mathrm{Cl}_2(k)$ is elementary in order to capitalize on genus theory. We characterize those $k$ for which $\operatorname{rank}\mathrm{Cl}_2(k^1) = 2$. We then show that the length of their 2-class field towers is 2. This raises a natural question.

> *Do there exist number fields $K$ for which* $\operatorname{rank}\mathrm{Cl}_p(K^1) = 2$, *but for which* $\mathrm{Cl}_p(K^2)$ *is not trivial?*

To carry out this program, we first refer to [1] for a classification of those imaginary quadratic $k$ with $\mathrm{Cl}_2(k) \simeq (2,2,2)$. (Here, $(2,2,2)$ denotes a group which is a direct sum of 3 cyclic groups of order 2.) Next, we use Koch's Theorem 1 of [17] to characterize those $k$ (as above) for which $\mathrm{Cl}_2(k^1)$ has rank 2, respectively 3, as a module over the integral group ring $\mathbb{Z}[\operatorname{Gal}(k^1/k)]$. In the cases where the module rank is 2, we sieve out further cases where $\operatorname{rank}\mathrm{Cl}_2(k^1) \geq 3$; and for the remaining cases, we use group theory to show that $\operatorname{rank}\mathrm{Cl}_2(k^1) = 2$ and then to glean enough information about the structure of $\operatorname{Gal}(k^2/k)$ to determine the length of the 2-class field tower of $k$. In many cases we give a presentation of $\operatorname{Gal}(k^\infty/k)$.

## 2. The Main Results

The following theorem classifies certain complex quadratic number fields for which the 2-class field tower terminates at the second step:

**Theorem 1.** *Consider the complex quadratic number fields $k$ with 2-class groups of type $(2,2,2)$ whose discriminants satisfy $d = d_1 d_2 d_3 d_4$, where $d_1$, $d_2$ and $d_3$ are negative prime discriminants. Let $G$ denote the Galois group of $k^2/k$, put $G_3 = [G, G']$.*

(1) *Among the unramified abelian extensions of degree 16 over $k(\sqrt{d_4})$, there is a unique extension $K$ for which $\mathrm{Cl}_2(K)$ has maximal rank; moreover, $K/k$ is normal and $\operatorname{Gal}(K/k) \simeq G/G_3$.*

(2) *The following assertions are equivalent:*
   i) $k^2 \neq k^3$;
   ii) $\operatorname{rank}\mathrm{Cl}_2(k^1) = 3$;
   iii) $G/G_3 \simeq 32.033$, *as classified in* [14];
   iv) *The discriminant* $\operatorname{disc} k = d_1 d_2 d_3 d_4$ *is a product of four prime discriminants* $d_1, d_2, d_3 < 0$, $d_4 > 0$ *such that* $(d_1/p_2) = (d_2/p_3) = (d_3/p_1) = (d_1/p_4) = -1$, $(d_4/p_2) = (d_4/p_3) = 1$, *where $p_i$ is the unique prime dividing $d_i$.*

   *If $k^2 = k^3$, on the other hand, then $\operatorname{rank}\mathrm{Cl}_2(k^1) = 2$.*

These results will be verified by going through the possible groups $G/G_3$. It turns out that the degree of difficulty is measured by the quantity $\rho$ defined as follows: first recall that a factorization $d \cdot d'$ of discriminants of quadratic fields is called a $C_4$-factorization of $dd'$ if $(d/p') = (d'/p) = +1$ for all primes $p \mid d$ and $p' \mid d'$. Let $\{d_1, d_2, d_3, d_4\}$ be the set of prime discriminants dividing $d = \operatorname{disc} k$;

TABLE 1

| $\rho$ | $G/G_3$ | $r_2(k^1)$ | $R_2$ | $r_4(k^1)$ | $\ell$ | $\mathcal{M}(G/G_3)$ |
|---|---|---|---|---|---|---|
| 0 | $32.040, 32.041$ | 2 | 3 | 0 | 2 | $(2)$ |
| 1 | $32.035, 32.037, 32.038$ | 2 | 3 | 1 | 2 | $(2,2)$ |
| 2 | $32.036$ | 2 | 3 | 2 | 2 | $(2,2,2)$ |
| 3 | $32.033$ | 3 | 4 | $2 \leq r_4 \leq 3$ | $3 \leq \ell \leq \infty$ | $(2,2,2,2)$ |

then $\rho$ is defined as the number of (independent) $C_4$-factorizations constructed out of these factors.

Table 1 reveals a close correlation between $\rho$, the 4-rank $r_4(k^1)$ of $\text{Cl}(k^1)$, and the rank of the Schur multiplier of $G/G_3$, as well as between the rank $r_2(k^1)$ of $\text{Cl}_2(k^1)$, the maximum $R_2$ of the 2-ranks of the class groups of the unramified quadratic extensions $k_j$ of $k$, and the length $\ell$ of the 2-class field tower of $k$. These unexpected results suggest that there are some structural relations between these invariants that are yet to be discovered.

## 3. Class Number Formulas and Unit indices

One of the major problems in applying class number formulas for multiquadratic fields is the computation of unit indices in fields of large degree. For a number field $L$, let $E_L$ denote its unit group and $W_L$ its subgroup of roots of unity. For CM-fields $L$ with maximal real subfields $K$, the Hasse unit index $Q(L) = (E_L : W_L E_K)$ can often be computed easily (see [21]). On the other hand, if $L/k$ is an elementary abelian 2-extension, the index $q(L/k) = (E_L : \prod e_i)$, where the product is over all quadratic subextensions $k_i/k$ of $L/k$ and where $e_i$ denotes the unit group of $k_i$, is in general much harder to compute. The following result is often helpful in reducing the necessary amount of computation:

**Proposition 1.** *Let $L$ be the compositum of quadratic number fields, and assume that $L$ does not contain a primitive eighth root of unity. If $L$ is complex and if $K$ denotes its maximal real subfield, then $q(L/\mathbb{Q}) = Q(L)q(K/\mathbb{Q})$.*

*Proof.* We have $q(L/\mathbb{Q}) = (E_L : \prod e_i) = (E_L : W_L E_K)(W_L E_K : \prod e_i)$, because the unit groups $e_i$ of real quadratic subfields are contained in $E_K$ and because $W_L = \prod e_j$, where the product is over the complex quadratic subfields $k_j$ of $L$ (here we have used that $L$ does not contain $\zeta_8$). Thus $q(L/\mathbb{Q}) = Q(L)(W_L E_K : \prod e_i) = (W_L E_K : W_L \prod e_r)$, where the last product is over all real quadratic subfields of $L$. Using the index formula $(UG : UH) = (G \cap U : H \cap U)(G : H)$ for subgroups $U, G, H$ of finite index in some abelian group and observing that $W_L \cap E_K = W_L \cap \prod e_r = \{-1, +1\}$ we get $(W_L E_K : W_L \prod e_r) = (E_K : \prod e_r) = q(K/\mathbb{Q})$. This proves our claim. □

This result will turn out to be useful in evaluating the class number formula for composita $L$ of quadratic number fields: if $(L : \mathbb{Q}) = 2^m$, then

$$(1) \qquad h(L) = 2^{-v} q(L/\mathbb{Q}) \prod h_i,$$

where the product is over the class numbers of all the quadratic fields in $L$ and where $v$ is defined by

$$v = \begin{cases} m(2^{m-1} - 1) & \text{if } L \text{ is real,} \\ (m-1)(2^{m-2} - 1) + 2^{m-1} - 1 & \text{if } L \text{ is complex.} \end{cases}$$

See e.g. Wada [28]. The following table gives the values of $v$ for the cases occurring in this paper:

| $L \backslash m$ | 2 | 3 | 4 |
|---:|:---:|:---:|:---:|
| real | 2 | 9 | 28 |
| complex | 1 | 5 | 16 |

**Ambiguous Class Number Formula.** We shall repeatedly make use of the ambiguous class number formula; let us recall the relevant notions and results here. For a cyclic extension $K/F$, $\text{Gal}(K/F) = \langle \sigma \rangle$ acts on $\text{Cl}(K)$ and thus defines the subgroup $\text{Am}(K/F)$ of $\text{Cl}(K)$ via the exact sequence

$$1 \longrightarrow \text{Am}(K/F) \longrightarrow \text{Cl}(K) \xrightarrow{1-\sigma} \text{Cl}(K)^{1-\sigma} \longrightarrow 1$$

as the subgroup of invariant (ambiguous) ideal classes. The following formula for $\#\text{Am}(K/F)$ is well known:

**Proposition 2.** *Let $K/F$ be a cyclic extension of prime degree $p$. Then*

$$\#\text{Am}(K/F) = h(F) \frac{p^{t-1}}{(E:H)},$$

*where $t$ is the number of (finite or infinite) primes of $F$ which are ramified in $K/F$, $E = E_F$ is the unit group of $F$, and $H = E \cap N_{K/F} K^\times$ is its subgroup of units which are norms of elements of $K$.*

For a proof, see [19].

We shall be interested in the case where $p = 2$ and $h(F)$ is odd. Let $\text{Am}_2$ denote the Sylow 2-subgroup of $\text{Am}$. Then it is known that $\#\text{Am}_2(K/F) = 2^e$ where $e$ equals the 2-rank of $\text{Cl}(K)$.

## 4. Imaginary Quadratic Fields $k$ with $\text{Cl}_2(k) \simeq (2,2,2)$.

In [1] a complete classification of imaginary quadratic number fields with $\text{Cl}_2(k) \simeq (2,2,2)$ is given. We now reduce the number of cases listed in [1] and represent each case by a graph. Let $d_k = \text{disc } k$ and assume $\text{Cl}_2(k) \simeq (2,2,2)$. By genus theory $d_k = d_1 d_2 d_3 d_4$ where $d_i$ are prime discriminants. Since $d_k < 0$, either exactly 1 or exactly 3 of the $d_i < 0$. For each of these two possibilities, we consider whether or not $d_k \equiv 4 \bmod 8$. From this we classify $d_k$ as follows:

**Type 1:** $d_k \not\equiv 4 \bmod 8$ and $d_i > 0$ $(i = 1, 2, 3)$, $d_4 < 0$;
**Type 2:** $d_k \not\equiv 4 \bmod 8$ and $d_i < 0$ $(i = 1, 2, 3)$, $d_4 > 0$;
**Type 3:** $d_k \equiv 4 \bmod 8$ and $d_i > 0$ $(i = 1, 2, 3)$, $d_4 = -4$;
**Type 4:** $d_k \equiv 4 \bmod 8$ and $d_i < 0$ $(i = 1, 2)$, $d_3 > 0$, $d_4 = -4$.

In order to determine when $\text{Cl}_2(k)$ is elementary, we use graphs on the primes dividing $\text{disc } k$ to do some convenient bookkeeping, cf. [6]. Let $p_i$ denote the prime dividing $d_i$ for $i = 1, 2, 3, 4$, i.e. $d_i = p_i^*$ where $p^* = (-1)^{(p-1)/2} p$ if $p$ is an odd prime and $2^* = 8, -4,$ or $-8$. Then we draw an arrow from $p_i$ to $p_j$ iff $(d_i/p_j) = -1$, where (./.) denotes the Kronecker symbol. If there is an arrow in both directions, we simply draw a line segment between the two primes.

The table at the end of the article gives a list of the discriminants of all imaginary quadratic fields with $\mathrm{Cl}_2(k) \simeq (2,2,2)$. The number and letter in the first column refer to the Type of $d_k$ given above and the subcases which correspond (as closely as possible) to the subcases in [1]. The second column gives the graph of the field. The four points in each graph represent the primes $p_1, p_2, p_3, p_4$, ordered as in the list of Types above, such that $p_1$ is the upper left point, $p_2$ upper right, $p_3$ lower left, and $p_4$ lower right. For example, 1A refers to a discriminant of the form $d_k = d_1 d_2 d_3 d_4 \not\equiv 4 \bmod 8$ with $d_i > 0$ for $i = 1,2,3$, and $d_4 < 0$, and such that $(d_1/p_4) = (d_2/p_4) = (d_3/p_4) = -1$, and $(d_i/p_j) = 1$ for $i,j = 1,2,3$, $i \neq j$.

## 5. The $\mathbb{Z}[\mathrm{Gal}(k^1/k)]$-rank of $\mathrm{Cl}_2(k^1)$

Let $K$ be a number field and $\Lambda$ the integral group ring $\mathbb{Z}[\mathrm{Gal}(K^1/K)]$. Then $\mathrm{Cl}_2(K^1)$ has a natural structure as a $\Lambda$-module induced by the action of the Galois group on the ideals of $K^1$.

Let $G = \mathrm{Gal}(K^2/K)$; let $G' = G_2 = [G,G]$ be the commutator subgroup of $G$; let $G'' = (G')'$; and finally let $G_n$ be defined inductively as $G_{n+1} = [G, G_n]$. Then by class field theory, we have $\mathrm{Cl}_2(K) \simeq \mathrm{Gal}(K^1/K) \simeq G/G'$, $\mathrm{Cl}_2(K^1) \simeq \mathrm{Gal}(K^2/K^1) = G'$. Moreover, the action of $\mathrm{Gal}(K^1/K)$ on $\mathrm{Cl}_2(K^1)$ corresponds, via the Artin map, to group conjugation of $G/G'$ on $G'$. Now suppose $G$ is a $p$-group; then by the work of Furtwängler [11] (along with Theorem 2.81 of [13]) it is known that the $\Lambda$-rank of $G'$ is equal to the $p$-rank of $G'/G_3$. This fact will be used in the proof of the following theorem.

**Theorem 2.** *Let $k$ be an imaginary quadratic number field such that $\mathrm{Cl}_2(k) \simeq (2,2,2)$.*
  *a) If $\mathrm{disc}\, k$ is divisible by three positive prime discriminants, then*
     $\mathrm{Cl}_2(k^1)$ *has rank 3 as a $\Lambda$-module.*
  *b) If $\mathrm{disc}\, k$ is divisible by one positive prime discriminant only, then*
     $\mathrm{Cl}_2(k^1)$ *has rank 2 as a $\Lambda$-module.*

*Proof.* (Sketch) Let $G = \mathrm{Gal}(k^2/k)$. Since the rank of $G/G_2$ is 3, the Burnside Basis Theorem implies that $G = \langle a_1, a_2, a_3 \rangle$ for some $a_i$ in $G$. Let $c_{ij} = [a_i, a_j] = a_i^{-1} a_j^{-1} a_i a_j$. Then $G' = \langle c_{12}, c_{13}, c_{23}, G_3 \rangle$, cf. [13]. Koch [17] has given a presentation of $G/H$, where $H$ is the subgroup $(G^2 G_2)^2 [G, G^2 G_2]$, in particular for imaginary quadratic fields. As the exponent of $G/G'$ is 2, $G_2^2 \subseteq G_3$ (for $[a,b]^2 \equiv [a^2, b] \equiv 1 \bmod G_3$). Thus Koch's presentation yields a presentation of $G/G_3$ which is given as follows (cf. [17] for details). For the prime discriminants $d_i = p_i^*$, let $d_k = d' d_i^\nu$ for $\nu = 0, 1$. Define $[d_k, p_i] \in \{0, 1\}$ by the relation

$$(-1)^{[d_k, p_i]} = (d'/p_i).$$

(Notice that $[d_i, p_i] = 0$.) Then for $k$ satisfying the assumption of the theorem, $G/G_3$ is presented as $\langle s_1, s_2, s_3 \rangle$ such that

$$\begin{aligned}
s_1^{2\delta_1} s_2^{2\nu_{21}} s_3^{2\nu_{31}} &= t_{12}^{\nu_{21}} t_{13}^{\nu_{31}} \\
s_1^{2\nu_{12}} s_2^{2\delta_2} s_3^{2\nu_{32}} &= t_{12}^{\nu_{12}} t_{23}^{\nu_{32}} \\
s_1^{2\nu_{13}} s_2^{2\nu_{23}} s_3^{2\delta_3} &= t_{13}^{\nu_{13}} t_{23}^{\nu_{23}} \\
s_1^{2\mu_1} s_2^{2\mu_2} s_3^{2\mu_3} &= 1 \\
G_3 &= 1,
\end{aligned}$$

where $\nu_{ij} = [d_i, p_j]$, $\delta_j = [d_k, p_j]$, $\mu_j = [d_j, p_4]$, and $t_{ij} = [s_i, s_j]$ the commutator of $s_i$ with $s_j$. The third column in the table at the end of this article gives a presentation of $G/G_3$ in each case. The fourth column indicates a number which represents the group classified in Hall and Senior [14] to which $G/G_3$ is isomorphic.

The table thus shows that $G'/G_3$ has rank 3 for $d_k$ of Types 1 and 3, and rank 2 for Types 2 and 4. □

By this theorem we see that if $\operatorname{rank} \operatorname{Cl}_2(k^1) = 2$, then the discriminant $d_k$ is of Type 2 or 4. In the next section, we describe more fully the structure of some of the groups $G$ listed in the table for $d_k$ of these two types. Notice by the table that in these two cases $G/G_3 = 32.033, 32.035, 32.036, 32.037, 32.038, 32.040, 32.041$.

## 6. Properties of Groups G with $G/G_3$ Isomorphic to $32.033, 32.035, 32.036, 32.037, 32.038, 32.040, 32.041$

We now consider a particular subclass of 2-groups $G$ such that $G/G' \simeq (2,2,2)$ with $G'/G_3$ of rank 2 (and hence $G'/G_3 \simeq (2,2)$, since $G_j/G_{j+1}$ is elementary). $G$ may be presented as $G = \langle a_1, a_2, a_3 \rangle$ where $G' = \langle c_{12}, c_{13}, G_3 \rangle$ with $c_{23} \in G_3$. (Recall that $c_{ij} = [a_i, a_j]$; and for that matter, let $c_{ijk} = [c_{ij}, a_k]$ and $c_{ijkl} = [c_{ijk}, a_l]$.)

We collect some general facts about commutators. First notice that by the *Witt identity* (see e.g. [8])
$$1 \equiv c_{ijk} c_{kij} c_{jki} \pmod{G''}$$
and so we have in particular
$$1 \equiv c_{123} c_{312} \pmod{G_4},$$
since $c_{231} \equiv 1 \pmod{G_4}$ as $c_{23} \in G_3$ and $G'' \subseteq G_4$. But from this and the observation that $\exp(G_j/G_{j+1}) = 2$, we see
$$c_{123} \equiv c_{312}^{-1} \equiv c_{132} \pmod{G_4}.$$
Similarly, notice that $c_{ijj} \equiv c_{jij} \pmod{G_4}$.

We now examine each of the seven groups listed above.

### 6.1. $G/G_3 \simeq 32.041$.

**Proposition 3.** *Let $G$ be a finite 2-group such that $G/G_3$ is isomorphic to group 32.041 as listed in Hall and Senior [14]. Then $G_3 = \langle 1 \rangle$ and so $G$ itself is isomorphic to group 32.041.*

*Proof.* In addition to the relations given above, we have
$$a_1^2 \equiv 1, a_2^2 \equiv c_{12} c_{13}, a_3^2 \equiv c_{12} \pmod{G_3},$$
(see [14]). But then,
$$1 = [a_2^2, a_2] \equiv [c_{12} c_{13}, a_2] \equiv c_{122} c_{132},$$
$$1 = [a_3^2, a_3] \equiv [c_{12}, a_3] \equiv c_{123},$$
$$1 \equiv [a_1^2, a_j] \equiv c_{1j}^2 c_{1j1},$$
$$c_{133} = [c_{13}, a_3] \equiv [a_2^2 a_3^2, a_3] \equiv c_{23}^2 \equiv 1,$$
$$c_{121} c_{131} \equiv [a_2^2, a_1] \equiv c_{12}^2 c_{122},$$
$$c_{121} \equiv [a_3^2, a_1] \equiv c_{13}^2 c_{131} \pmod{G_4}.$$

This implies that
$$c_{ijk} \equiv 1 \pmod{G_4}.$$
By Theorem 2.81 of [13],
$$G_3 = \langle c_{121}, c_{131}, c_{122}, c_{132}, c_{123}, c_{133}, G_4 \rangle = G_4.$$
Thus $G_3 = \langle 1 \rangle$, since $G$ is nilpotent. This establishes the proposition. □

Using the table at the end and the above proposition, we can characterize those imaginary quadratic fields $k$ for which $\mathrm{Gal}(k^\infty/k)$ is isomorphic to group 32.041.

**Proposition 4.** *Let $k$ be an imaginary quadratic number field with discriminant $d_k$. Then $\mathrm{Gal}(k^\infty/k) \simeq 32.041$ if and only if there is a factorization of $d_k = d_1 d_2 d_3 d_4$ into distinct prime discriminants $d_j$, divisible by the unique primes $p_j$, satisfying*
  (1) $d_k \not\equiv 4 \bmod 8$,
  (2) $d_i < 0$, *for* $i = 1, 2, 3$ *and* $d_4 > 0$,
  (3) $(d_i/p_4) = -1$ *for* $i = 1, 2, 3$ *and* $(d_1/p_2) = (d_2/p_3) = (d_3/p_1) = -1$.

## 6.2. $G/G_3 \simeq 32.040$.

**Proposition 5.** *Let $G$ be a finite 2-group such that $G/G_3$ is isomorphic to group 32.040. Then $G_3 = \langle 1 \rangle$.*

*Proof.* $G$ has the additional relations
$$a_1^2 \equiv c_{12}, a_2^2 \equiv c_{12}c_{13}, a_3^2 \equiv c_{13} \pmod{G_3}.$$
Consequently,
$$1 = [a_1^2, a_1] \equiv c_{121},$$
$$1 = [a_3^2, a_3] \equiv c_{133},$$
$$1 = [a_2^2, a_2] \equiv c_{122}c_{132},$$
$$c_{123} = [c_{12}, a_3] \equiv [a_2^2 a_3^2, a_3] \equiv c_{23}^2 \equiv 1,$$
$$c_{122} \equiv [a_1^2, a_2] \equiv c_{12}^2 c_{121},$$
$$c_{121}c_{131} \equiv [a_2^2, a_1] \equiv c_{12}^2 c_{122},$$
$$c_{131} \equiv [a_3^2, a_1] \equiv c_{13}^2 c_{133} \pmod{G_4}.$$
All this yields,
$$c_{ijk} \equiv 1 \pmod{G_4}.$$
Hence $G_3 = \langle 1 \rangle$. □

We now characterize those imaginary quadratic fields $k$ for which $\mathrm{Gal}(k^\infty/k)$ is isomorphic to group 32.040.

**Proposition 6.** *Let $k$ be an imaginary quadratic number field with discriminant $d_k$. Then $\mathrm{Gal}(k^\infty/k) \simeq 32.040$ if and only if there is a factorization of $d_k = -4qq'p$, with $p, q, q'$ distinct primes satisfying*
  (1) $q \equiv q' \equiv 3 \bmod 4$ *and* $p \equiv 1 \bmod 4$,
  (2) $(q/p) = (q'/p) = -1$, $(q'/q) = +1$, $q \equiv 3 \bmod 8$, $q' \equiv 7 \bmod 8$, *and* $p \equiv 5 \bmod 8$.

6.3. $\mathbf{G/G_3 \simeq 32.038}$.

**Proposition 7.** *Let $G$ be a finite 2-group such that $G/G_3$ is isomorphic to group 32.038. Then $G' \simeq (2, 2^m)$ for some $m \geq 1$.*

*Proof.* The additional relations are
$$a_1^2 \equiv a_3^2 \equiv 1, \ a_2^2 \equiv c_{13} \pmod{G_3},$$
whence
$$1 \equiv [a_1^2, a_j] \equiv c_{1j}^2 c_{1j1},$$
$$1 \equiv [a_3^2, a_j] \equiv c_{3j}^2 c_{3j3},$$
$$1 = [a_2^2, a_2] \equiv c_{132},$$
$$c_{133} \equiv [a_2^2, a_3] \equiv c_{23}^2 \equiv 1,$$
$$c_{131} \equiv [a_2^2, a_1] \equiv c_{12}^2 c_{122} \pmod{G_4}.$$

But then
$$c_{132} \equiv c_{123} \equiv c_{133} \equiv c_{131} \equiv c_{13}^2 \equiv 1, \ c_{122} \equiv c_{121} \equiv c_{12}^2 \pmod{G_4},$$
which implies (as above) that $G_3 = \langle c_{12}^2, G_4 \rangle$ and more generally
$$G_j = \langle c_{12}^{2^{j-2}}, G_{j+1} \rangle.$$
But since $G$ is nilpotent $G_j$ is trivial for $j$ large enough. This implies that $G_3 = \langle c_{12}^2 \rangle$. But then $G'/G'' \simeq (2, 2^m)$. By Theorem 1 of [7], we see that $G''$ is trivial. Thus the proposition follows. □

We now wish to compute $G = \mathrm{Gal}(k^\infty/k)$ for those $k$ with $G/G_3 \simeq 32.038$. To this end, we start with a lemma.

**Lemma 1.** *Let $G$ be a finite 2-group such that $G/G_3 \simeq 32.038$. Hence, we know $G = \langle a_1, a_2, a_3 \rangle$ with*
$$a_1^2 \equiv a_3^2 \equiv c_{23} \equiv 1, \ a_2^2 \equiv c_{13} \mod G_3.$$
*Let $A = \langle a_2, a_3, G' \rangle$, and $B = \langle a_3, a_1 a_2, G' \rangle$. Finally, suppose that $\ker t_B = \langle \overline{a_1 a_2}, \overline{a_1 a_3} \rangle$ or $\langle \overline{a_2}, \overline{a_1 a_3} \rangle$, where $t_H$ represents the transfer map from $G$ to a subgroup $H$, and $\overline{a} = aG'$.*

*Then $G \simeq \Gamma_n^{(38)}$, for some integer $n$ with $n \geq 2$, where*
$$\Gamma_n^{(38)} = \langle a_1, a_2, a_3 : \ a_1^4 = a_2^4 = a_3^2 = 1, \ a_1^2 = c_{12}^{2^{n-1}}, \ a_2^2 = c_{13}, \ c_{23} = a_1^2 \ \rangle.$$

*Proof.* First notice that if we replace $a_2$ by $a_2 a_3$ we may assume without loss of generality that
$$\ker t_B = \langle \overline{a_1 a_2}, \overline{a_1 a_3} \rangle.$$
Next notice that
$$A' = G_3, \ B' = \langle c_{13} c_{23} c_{132}, c_{123}, c_{133}, c_{121} c_{122} c_{1212}, c_{131} c_{132} c_{1312}, B_3 \rangle.$$
We now list the values of the transfer maps $t_H$, for $H = A, B$. First recall that if $(G : H) = 2$ so that $G = H \cup Hz$, then $t_H$ is the homomorphism from $G/G'$ into $H/H'$ determined by
$$t_H(\overline{h}) = hz^{-1}hzH' = h^2[h, z]H', \qquad t_H(\overline{hz}) = (hz)^2 H'$$
for any $h \in H$.

| $i$ | $t_A(\overline{a_i})$ | $t_B(\overline{a_i})$ |
|---|---|---|
| 1 | $a_1^2 A' = G_3$ | $a_1^2 B'$ |
| 2 | $a_2^2 c_{12}^{-1} A' = c_{12} c_{13} G_3$ | $a_2^2 B'$ |
| 3 | $a_3^2 c_{13}^{-1} A' = c_{13} G_3$ | $a_3^2 c_{13}^{-1} B'$ |

Notice that
$$\ker t_A = \langle \overline{a_1} \rangle.$$
By our assumptions, $t_B(\overline{a_1}) \neq B'$. We claim that this implies that $a_1^4 = 1$ and $B' \cap G_3 = \langle 1 \rangle$. To see all this, notice by the table and assumption that $a_1^2 \notin B'$. On the other hand, $B' = t_B(\overline{1}) = t_B(\overline{a_1}^2) = a_1^4 B'$. Thus $a_1^4 \in B' \cap G_3$. Since $G_3 = \langle c_{12}^2 \rangle$, we see $B' \cap G_3 = \langle a_1^4 \rangle$. We claimed $a_1^4 = 1$. Suppose not. Then let $a_1^2 = c_{12}^{2^k x}$ for some $k \geq 1$ and $x$ odd (recall that $a_1^2 \in G_3$). Hence, since $G_j = \langle c_{12}^{2^{j-2}} \rangle$ for $j > 2$, $a_1^2 \in G_{k+2} \setminus G_{k+3}$ and thus $a_1^4 \in G_{k+3} \setminus G_{k+4}$ (since $a_1^4 \neq 1$). Now notice that
$$1 = [a_1^2, a_1] = [c_{12}^{2^k x}, a_1] = c_{121}^{2^k x}, \quad c_{12}^2 c_{121} = [a_1^2, a_2] = [c_{12}^{2^k x}, a_2] = c_{122}^{2^k x},$$
implying that $c_{12}^2 = c_{122}^{2^k x} c_{121}^{-1}$ and
$$a_1^4 = (c_{12}^2)^{2^k x} = c_{122}^{2^{2k} x^2} c_{121}^{-2^k x} = c_{122}^{2^{2k} x^2} \in G_{2k+3} \subseteq G_{k+4},$$
a contradiction. Therefore, $a_1^4 = 1$, as claimed, and so $B' \cap G_3$ is trivial. Also notice that since $a_1^2 \in G_3 = \langle c_{12}^2 \rangle$ and $a_1^2 \neq 1$, $a_1^2 = c_{12}^{2^{n-1}}$, where $2^n$ is the order of $c_{12}$.

From the argument above we have $c_{12}^2 = c_{121}^{-1}$ and thus $c_{1212} = [c_{12}^{-2}, a_2] = c_{122}^{-2}$. Moreover, since $B' \cap G_3 = \langle 1 \rangle$, we get
$$c_{123} = c_{133} = c_{121} c_{122} c_{1212} = c_{131} c_{133} c_{1312} = 1,$$
and so $c_{122} = c_{12}^{-2}$.

Next, we have $c_{233} = 1$, for $c_{23} = c_{12}^{2x}$ (since $c_{23} \in G_3$) and thus $c_{233} = c_{123}^{2x} = 1$, since $c_{123} = 1$.

Notice $a_2^4 = 1$, for $B' = t_B(\overline{a_2}^2) = a_2^4 B'$ and so $a_2^4 \in B' \cap G_3 = \langle 1 \rangle$.

Next we claim $c_{132} = c_{13}^{-2}$, for since $a_2^2 \equiv c_{13} \mod G_3$ and $G_3 = \langle c_{12}^2 \rangle$, we have $a_2^2 = c_{13} c_{12}^{2x}$ and so $1 = [a_2^2, a_2] = c_{132} c_{122}^{2x} = c_{132} c_{12}^{-4x}$. On the other hand, $1 = a_2^4 = c_{13}^2 c_{12}^{4x}$, and so cancellation yields claim.

This all implies $B' = \langle c_{13}^{-1} c_{23} \rangle$ and since $(c_{13}^{-1} c_{23})^2 \in B' \cap G_3$ we have $c_{13}^2 = c_{23}^2$.

Now we claim $a_1^2 a_2^2 = c_{13} c_{23}^{-1}$ and $a_1^2 a_3^2 = c_{23}$; for $B' = t_B(\overline{a_1 a_2}) = a_1^2 a_2^2 B'$ implying that $a_1^2 a_2^2 \in B' \setminus G_3 = \{c_{13} c_{23}^{-1}\}$ which gives the first equality. For the second, $B' = t_B(\overline{a_1 a_3}) = a_1^2 a_3^2 c_{23}^{-1} B'$, whence $a_1^2 a_3^2 c_{23}^{-1} \in B' \cap G_3$, implying the second part of the claim.

Summarizing, we have
$$c_{12}^{2^n} = a_1^4 = a_2^4 = 1, \ a_1^2 = c_{12}^{2^{n-1}}, \ a_2^2 = a_1^2 c_{13} c_{23}^{-1}, \ a_3^2 = a_1^2 c_{23}, \ c_{13}^4 = 1.$$

There are two cases to consider:

*Case 1* Suppose $c_{13}^2 = 1$.
Then $c_{23}^2 = 1$. If $c_{23} = c_{12}^{2^{n-1}}$, then we obtain the desired presentation. If $c_{23} = 1$, then replacing $a_3$ by $a_3 c_{12}^{2^{n-2}}$ gives the same presentation.

*Case 2* Suppose $c_{13}^2 \neq 1$.
Let $c_{23} = c_{12}^{2x}$. Then $c_{13}^2 = c_{23}^2 = c_{12}^{4x}$. If we replace $a_3$ by $a_3 c_{12}^{-x}$, then we are reduced to *Case 1*.

This establishes the lemma. □

**Proposition 8.** *Let $k$ be an imaginary quadratic number field with discriminant $d_k$. Furthermore, let $G = \mathrm{Gal}(k^\infty/k)$. Then $G/G_3 \simeq 32.038$ if and only if there is a factorization of $d_k = -4qq'p$, with $p, q, q'$ distinct primes satisfying*

(1) *$q \equiv q' \equiv 3 \bmod 4$ and $p \equiv 1 \bmod 4$,*
(2) *$(q/p) = +1$, $(q'/p) = -1$, $(q'/q) = +1$, $q \equiv 3 \bmod 8$, and $p \equiv 5 \bmod 8$.*

*Moreover, $G \simeq \Gamma_n^{(38)}$ where $n$ is determined by $2^n = h_2(-qp)$. In particular, $\mathrm{Cl}_2(k^1) \simeq (2, 2^n)$.*

*Proof.* The first part of the proposition follows immediately from the table at the end of this article. For proving $2^n = h_2(-qp)$, we compute a few class numbers and capitulation kernels:

| $j$ | $k_j$ | $\kappa_j$ | $h(k_j)$ | rank $C_j$ | $\mathrm{Cl}_2(k_j)$ | $NC_j$ | |
|---|---|---|---|---|---|---|---|
| 1 | $k(\sqrt{p})$ | $\langle[\mathfrak{p}]\rangle$ | 16 | 3 | $(2,2,4)$ | $\langle[2\mathfrak{p}],[\mathfrak{q}]\rangle$ | |
| 2 | $k(\sqrt{q'})$ | $\langle[2],[\mathfrak{pq}]\rangle$ | $2^{n+2}$ | 2 | $(2, 2^{n+1})$ | $\langle[2],[\mathfrak{q}]\rangle$ | |
| 3 | $k(\sqrt{-q})$ | $\langle[2\mathfrak{p}],[\mathfrak{q}]\rangle$ | 8 | 3 | $(2,2,2)$ | $\langle[\mathfrak{p}],[\mathfrak{q}]\rangle$ | |
| 4 | $k(\sqrt{-1})$ | $\langle[2],[\mathfrak{q}]\rangle$ | 8 | 2 | $(2,4)$ | $\langle[2\mathfrak{q}],[\mathfrak{p}]\rangle$ | $\langle[2],[\mathfrak{p}]\rangle$ |
| 5 | $k(\sqrt{-p})$ | $\langle[\mathfrak{p}],[\mathfrak{q}]\rangle$ | 8 | 2 | $(2,4)$ | $\langle[2],[\mathfrak{pq}]\rangle$ | $\langle[2\mathfrak{p}],[2\mathfrak{q}]\rangle$ |
| 6 | $k(\sqrt{q})$ | $\langle[2],[\mathfrak{q}]\rangle$ | 8 | 2 | $(2,4)$ | $\langle[2],[\mathfrak{p}]\rangle$ | $\langle[2\mathfrak{q}],[\mathfrak{p}]\rangle$ |
| 7 | $k(\sqrt{-q'})$ | $\langle[\mathfrak{p}],[\mathfrak{q}]\rangle$ | 8 | 2 | $(2,4)$ | $\langle[2\mathfrak{p}],[2\mathfrak{q}]\rangle$ | $\langle[2],[\mathfrak{pq}]\rangle$ |

Here $\kappa_j$ is the capitulation kernel in $k_j/k$, $C_j = \mathrm{Cl}_2(k_j)$, $NC_j = N_{k_j/k}\mathrm{Cl}_2(k_j)$, and $2^n$ is the 2-class number of $\mathbb{Q}(\sqrt{-pq})$. The left column in $NC_j$ denotes the case $q' \equiv 3 \bmod 8$, the right one $q' \equiv 7 \bmod 8$.

The group $G/G_3 \simeq 32.038$ has seven maximal subgroups; one of them has abelianization $(2,2,4)$ (thus $G/G_3 \simeq \mathrm{Gal}(K^1/k)$ for $K = k_1$), another has abelianization $(2, 2^{n+1})$, and the others have abelianizations of order 8.

Most of the computations involved in verifying this table are straight forward and left to the reader. Let us check that rank $\mathrm{Cl}_2(k_4) = 2$. This is done via the ambiguous class number formula. We distinguish two cases:

a) $q' \equiv 3 \bmod 8$: then 2 is inert in $F = \mathbb{Q}(\sqrt{-q'})$, hence $t = 3$, and $-1$ is a quadratic residue modulo $p\mathcal{O}_F$ and $q\mathcal{O}_F$. Thus $E = H$ and rank $\mathrm{Cl}_2(k_4) = 2$.

b) $q' \equiv 7 \bmod 8$: Here $2\mathcal{O}_F = \mathfrak{z}_1\mathfrak{z}_2$, hence $t = 4$, and we want to show that $-1$ is not a local norm at $\mathfrak{z}_1$. To this end observe that $\mathbb{Q}_2(\sqrt{-q'}) = \mathbb{Q}_2$ and $\mathbb{Q}_2(\sqrt{pq}) = \mathbb{Q}_2(\sqrt{-1})$; thus the completion of $k_4$ at $\mathfrak{z}_1$ is $\mathbb{Q}_2(\sqrt{-1})$, and $-1$ is not a norm in $\mathbb{Q}_2(\sqrt{-1})/\mathbb{Q}_2$. Therefore $(E : H) = 2$, and our claim follows.

The other entries in our table are checked similarly, with the exception of the claim that $\mathrm{Cl}_2(k_2) \simeq (2, 2^{n+1})$ which requires more care. In fact, consider $F =$ $\mathbb{Q}(\sqrt{-pq})$ and let $\mathrm{Cl}_2(F)$ be generated by the class of the ideal $\mathfrak{a}$, say. Then $\mathfrak{a}^{2^{n-1}} \sim \mathfrak{q}$, and since $\mathfrak{q}$ does not capitulate in $k_2/F$, we see that $\mathfrak{a}$ still has order $2^n$ in $\mathrm{Cl}_2(k_2)$. Since $k_2/F$ is ramified, class field theory guarantees the existence of an ideal $\mathfrak{A}$ in $k_2$ such that $N_{k_2/F}\mathfrak{A} \sim \mathfrak{a}$; letting $\sigma$ denote the nontrivial automorphism of $k_2/F$, this means that $\mathfrak{A}^{1+\sigma} \sim \mathfrak{a}$ in $\mathrm{Cl}_2(k_2)$. On the other hand, since $\mathrm{Cl}_2(k)$ has exponent 2, we deduce $\mathfrak{A}^{1+\tau} = \mathfrak{b}$ with $\mathfrak{b}^2 \sim 1$, where $\tau$ is the nontrivial automorphism of $k_2/k$. Finally, $\mathfrak{A}^{1+\sigma\tau} \sim 1$ since the class number of $\mathbb{Q}(\sqrt{q'})$ is odd. Multiplying these relations gives $\mathfrak{A}^2 \sim \mathfrak{A}^{1+\sigma}\mathfrak{A}^{1+\tau}\mathfrak{A}^{1+\sigma\tau} \sim \mathfrak{ab}$, and since $[\mathfrak{ab}]$ has order $2^n$, the ideal class of $\mathfrak{A}$ must have order $2^{n+1}$.

From the information on the capitulation kernels and orders of the 2-class groups, we see that $\mathrm{Gal}(k^2/k_1) = A$ and $\mathrm{Gal}(k^2/k_2) = B$ with $A$ and $B$ as given in Lemma 1. Notice that $[\mathfrak{p}]$ and $[\mathfrak{q}]$ are then identified with $\overline{a_1}$ and $\overline{a_3}$, respectively, under the identification of $\mathrm{Cl}_2(k)$ with $G/G'$. Consequently, $G$ satisfies the hypotheses of the lemma above and thus $G \simeq \Gamma_n^{(38)}$.

Finally, we need to show that $\#G' = 2h_2(-qp)$. To this end, we know by Proposition 7 that $\#G' = 2^{n+1}$ for some $n$. On the other hand, by the table and class number formula for biquadratic fields,
$$4h_2(-qp) = h_2(k_2) = (B:B') = 2^{n+2},$$
as desired.

This establishes the proposition. □

Numerical examples:

| $d$ | $p$ | $q$ | $q'$ | $\mathrm{Cl}_2(k^1)$ |
|---|---|---|---|---|
| $-660$ | 5 | 11 | 3 | $(4,2)$ |
| $-1092$ | 13 | 3 | 7 | $(4,2)$ |

### 6.4. $G/G_3 \simeq 32.037$.

**Proposition 9.** *Let $G$ be a finite 2-group such that $G/G_3$ is isomorphic to group 32.037. Then $G' \simeq (2, 2^m)$ for some $m \geq 1$.*

*Proof.* The additional relations are
$$a_1^2 \equiv c_{13}, a_2^2 \equiv 1, a_3^2 \equiv c_{13} \pmod{G_3}$$
and so
$$1 \equiv [a_2^2, a_1] \equiv c_{12}^2 c_{122},$$
$$1 = [a_1^2, a_1] \equiv [c_{13}, a_1] \equiv c_{131},$$
$$1 = [a_3^2, a_3] \equiv [c_{13}, a_3] \equiv c_{133},$$
$$c_{132} \equiv [a_3^2, a_2] \equiv c_{23}^2 \equiv 1,$$
$$c_{132} \equiv [a_1^2, a_2] \equiv c_{12}^2 c_{121} \pmod{G_4}.$$
Hence,
$$c_{121} \equiv c_{122} \equiv c_{12}^2, \quad c_{131} \equiv c_{133} \equiv c_{132} \equiv c_{123} \equiv 1 \pmod{G_4},$$
whence $G_3 = \langle c_{12}^2, G_4 \rangle$. The arguments as above then show that $G_3 = \langle c_{12}^2 \rangle$. Therefore, $G'$ is abelian with $G' \simeq (2, 2^m)$, as desired. □

**Lemma 2.** *Let $G$ be a finite 2-group such that $G/G_3 \simeq 32.037$; hence we have $G = \langle a_1, a_2, a_3 \rangle$ with*

$$a_1^2 \equiv a_3^2 \equiv c_{13}, \quad a_2^2 \equiv c_{23} \equiv 1 \bmod G_3.$$

*Let $B = \langle a_1 a_2, a_3, c_{12}, c_{13} \rangle$. Define*

$$\Gamma_{n,\varepsilon}^{(37)} = \langle a_1, a_2, a_3 : a_1^4 = a_2^4 = a_3^4 = c_{12}^{2^n} = 1, \ a_1^2 = c_{13}, \ a_3^2 = c_{13} a_2^{2\varepsilon}, \ a_2^2 = c_{12}^{2^{n-1}} = c_{23} \rangle,$$

*where $\varepsilon = 0, 1$.*

  (1) *If $\ker t_B = B/G'$, then $G \simeq \Gamma_{n,0}^{(37)}$ for some $n \geq 2$.*
  (2) *If $\ker t_B = \langle \overline{a_1 a_2}, \overline{a_1 a_3} \rangle$, then $G \simeq \Gamma_{n,1}^{(37)}$ for some $n \geq 2$.*

*Proof.* First, recall that $G_3 = \langle c_{12}^2 \rangle = \langle c_{122} \rangle$ and that $a_2^2 \in G_3$. Then in both cases of the lemma, $\overline{a_2} \notin \ker t_B$ and so the same argument as in the proof of Lemma 1 shows that $a_2^4 = 1$, $a_2^2 = c_{12}^{2^{n-1}}$, where $2^n$ is the order of $c_{12}$, and that $B' \cap G_3 = \langle 1 \rangle$. Hence, $B' = \langle c_{13} c_{23} c_{132} \rangle$ and $(c_{13} c_{23} c_{132})^2 = 1$ being in $B' \cap G_3$. For that matter

$$c_{123} = c_{133} = c_{121} c_{122} c_{1212} = c_{131} c_{132} c_{1312} = c_{233} = 1.$$

By assumption, $B' = t_B(\overline{a_1 a_2}) = a_1^2 a_2^2 B'$, which implies that $a_1^2 a_2^2 = c_{13} c_{23} c_{132}$, being an element in $B'$ but not $G_3$. Hence $a_1^4 = 1$.

If $B' = t_B(\overline{a_3}) = a_3^2 c_{13}^{-1} B'$, then $a_3^2 c_{13}^{-1} = 1$, being in $B' \cap G_3$. Hence $a_3^2 = c_{13}$ in this case. On the other hand, if $B' \neq t_B(\overline{a_3}) = a_3^2 c_{13}^{-1} B'$, then $a_3^2 c_{13}^{-1} \in G_3 \setminus B'$, whence $a_3^2 c_{13}^{-1} = c_{12}^{2x} \neq 1$. But since $B' = t_B(\overline{a_3^2})$, $a_3^4 c_{13}^{-2} = 1$, being in $B' \cap G_3$. Thus $a_3^2 = c_{13} c_{12}^{2^{n-1}}$ in this case.

Next, notice that in both cases of the lemma $c_{132} = [c_{13}, a_2] = [a_3^2, a_2] = c_{23}^{-2}$. Thus $c_{13} c_{23} c_{132} = c_{13} c_{23}^{-1}$, which in turn implies that $(c_{13} c_{23}^{-1})^2 = 1$ and so $c_{13}^2 = c_{23}^2$. Also notice that $c_{131}^{-1} = [a_1, c_{13}] = [a_1, a_3^2] = c_{13}^2 c_{133} = c_{13}^2$.

Summarizing, we have so far:

$$c_{12}^{2^n} = a_1^4 = a_2^4 = 1, \ a_2^2 = c_{12}^{2^{n-1}}, \ a_1^2 = a_2^2 c_{13} c_{23}^{-1}, \ a_3^2 = c_{13} a_2^{2\varepsilon}, \ c_{13}^4 = 1,$$

where $\varepsilon = 0$ if $t_B(\overline{a_3}) = B'$, and $\varepsilon = 1$ if not.

*Case 1.* Suppose $c_{13}^2 = 1$. If $c_{23} = c_{12}^{2^{n-1}}$, then we get the desired presentations. If $c_{23} = 1$, then replacing $a_3$ by $a_3 c_{12}^{2^{n-2}}$ yields the same presentations.

*Case 2.* Suppose $c_{13}^2 \neq 1$. If $c_{23} = c_{12}^{2x}$, then replacing $a_3$ by $a_3 c_{12}^{-x}$ yields Case 1. This establishes the lemma. □

**Proposition 10.** *Let $k$ be an imaginary quadratic number field with discriminant $d_k$. Furthermore, let $G = \mathrm{Gal}(k^\infty/k)$. Then $G/G_3 \simeq 32.037$ if and only if there is a factorization of $d_k = d_1 d_2 d_3 d_4$ into distinct prime discriminants satisfying*

  (1) $d_i < 0$, $(i = 1, 2, 3)$, *and* $d_4 > 0$,
  (2) $(d_1/p_4) = (d_2/p_4) = -1$, $(d_2/p_3) = (d_4/p_3) = 1$, $(d_3/p_2) = (d_2/p_1) = (d_1/p_3) = -1$
    *where $p_i$ is the unique prime dividing $d_i$.*

*Moreover, $G \simeq \Gamma_{n,0}^{(37)}$ if $d_k \equiv 4 \bmod 8$, and $G \simeq \Gamma_{n,1}^{(37)}$ if $d_k \not\equiv 4 \bmod 8$ where $n$ is determined by $2^n = h_2(d_3 d_4)$.*

*Proof.* As usual the congruence conditions come from our table at the end of the article. Now, let $\mathfrak{p}_j$ be the unique prime ideal of $k$ dividing $d_j$, for $j = 1, 2, 3, 4$. Consider the fields $k_4$ and $k_5$ in the table below. If we let $G$ be presented as in Proposition 9 and let $A = \langle a_2, a_3, G' \rangle$ and $B = \langle a_1a_2, a_3, G' \rangle$, then from this we have that

$$\operatorname{Gal}(k^\infty/k_4) = A, \quad \operatorname{Gal}(k^\infty/k_5) = B.$$

Assume first that $d_k \equiv 4 \bmod 8$. Then the congruence conditions imply that $d_3 = -4$, whence $\mathfrak{p}_1\mathfrak{p}_2\mathfrak{p}_4 \sim 1$, (where the equivalence is in $k$). But then $\kappa_5 = NC_5$, where $NC_j = N_{k_j/k}(\operatorname{Cl}_2(k_j))$, and so we have

$$\ker t_B = B/G'.$$

Hence by Lemma 2, we see that

$$\operatorname{Gal}(k^\infty/k) \simeq \Gamma_{n,0}^{(37)}.$$

Now assume that $d_k \not\equiv 4 \bmod 8$. At this point we need to determine $\ker t_B$ in terms of $a_1, a_2, a_3$. Since $A' = G_3$ as is easily seen, we see that $\ker t_A = \langle \overline{a_3} \rangle$ (see the proof of Lemma 1). Hence $\overline{a_3}$ must correspond to $[\mathfrak{p}_4]$ under the Artin map. Thus $\overline{a_3} \notin \ker t_B$. On the other hand, $[\mathfrak{p}_2] \in \kappa_5 \cap NC_5$, implying that $\overline{a_1a_2}$ or $\overline{a_1a_2a_3}$ in $\ker t_B$. If we replace $a_1$ by $a_1a_3$, if necessary, we may assume $\overline{a_1a_2} \in \ker t_B$. We then have the following two possibilities: $\ker t_B = \langle \overline{a_1a_2}, \overline{a_1a_3} \rangle$ or $\ker t_B = \langle \overline{a_1}, \overline{a_2} \rangle$. We claim that the first must occur. For, $\overline{a_2} \in \ker t_D$, for any maximal subgroup $D$ of $G$ other than $A$ and $B$ (exercise to the reader). Since $\langle [\mathfrak{p}_3] \rangle = \kappa_3 \cap \kappa_6$, we see that $[\mathfrak{p}_3]$ corresponds to $\overline{a_2}$. Since $[\mathfrak{p}_3] \notin \kappa_5$, the claim follows. Thus by Lemma 2, we have $G \simeq \Gamma_{n,1}^{(37)}$.

The rest of the proposition follows as in the proof of Proposition 8. This establishes the proposition. □

| $j$ | $k_j$ | $\kappa_j$ | $h(k_j)$ | rank $C_j$ | $\operatorname{Cl}_2(k_j)$ | $NC_j$ |
|---|---|---|---|---|---|---|
| 1 | $k(\sqrt{d_1})$ | $\langle [\mathfrak{p}_1], [\mathfrak{p}_3] \rangle$ | 8 | 2 | $(2,4)$ | $\langle [\mathfrak{p}_2], [\mathfrak{p}_3\mathfrak{p}_4] \rangle$ |
| 2 | $k(\sqrt{d_2})$ | $\langle [\mathfrak{p}_2], [\mathfrak{p}_3] \rangle$ | 8 | 2 | $(2,4)$ | $\langle [\mathfrak{p}_2], [\mathfrak{p}_3] \rangle$ |
| 3 | $k(\sqrt{d_3})$ | $\langle [\mathfrak{p}_1\mathfrak{p}_2], [\mathfrak{p}_3] \rangle$ | 8 | 3 | $(2,2,2)$ | $\langle [\mathfrak{p}_2\mathfrak{p}_3], [\mathfrak{p}_4] \rangle$ |
| 4 | $k(\sqrt{d_4})$ | $\langle [\mathfrak{p}_4] \rangle$ | 16 | 3 | $(2,2,4)$ | $\langle [\mathfrak{p}_1\mathfrak{p}_2], [\mathfrak{p}_3] \rangle$ |
| 5 | $k(\sqrt{d_1d_2})$ | $\langle [\mathfrak{p}_1], [\mathfrak{p}_2] \rangle$ | $2^{n+2}$ | 2 | $(2, 2^{n+1})$ | $\langle [\mathfrak{p}_2], [\mathfrak{p}_4] \rangle$ |
| 6 | $k(\sqrt{d_1d_3})$ | $\langle [\mathfrak{p}_1], [\mathfrak{p}_3] \rangle$ | 8 | 2 | $(2,4)$ | $\langle [\mathfrak{p}_1], [\mathfrak{p}_3] \rangle$ |
| 7 | $k(\sqrt{d_1d_4})$ | $\langle [\mathfrak{p}_2], [\mathfrak{p}_3] \rangle$ | 8 | 2 | $(2,4)$ | $\langle [\mathfrak{p}_2\mathfrak{p}_3], [\mathfrak{p}_3\mathfrak{p}_4] \rangle$ |

In the table above, we have the relations

$$[\mathfrak{p}_1\mathfrak{p}_2\mathfrak{p}_3] = [\mathfrak{p}_4] \quad \text{if} \quad d \not\equiv 4 \bmod 8,$$
$$[\mathfrak{p}_1\mathfrak{p}_2] = [\mathfrak{p}_4] \quad \text{if} \quad d \equiv 4 \bmod 8.$$

Numerical Examples

| $d$ | $d_1$ | $d_2$ | $d_3$ | $d_4$ | $h_2(d_3d_4)$ |
|---|---|---|---|---|---|
| $-1155$ | $-3$ | $-7$ | $-11$ | 5 | 4 |
| $-1428$ | $-3$ | $-7$ | $-4$ | 17 | 4 |
| $-3003$ | $-7$ | $-11$ | $-3$ | 13 | 4 |
| $-3444$ | $-3$ | $-7$ | $-4$ | 41 | 8 |

6.5. **$G/G_3 \simeq 32.035$.**

For the moment we skip over group 32.036 and consider the simpler case 32.035.

**Proposition 11.** *Let $G$ be a finite $2$-group such that $G/G_3$ is isomorphic to group 32.035. Then $G' \simeq (2, 2^m)$ for some $m \geq 1$.*

*Proof.* The additional relations are
$$a_1^2 \equiv c_{13}, a_2^2 \equiv c_{12}, a_3^2 \equiv c_{13} \pmod{G_3}$$
and so
$$1 = [a_1^2, a_1] \equiv c_{131},$$
$$1 = [a_2^2, a_2] \equiv [c_{12}, a_2] \equiv c_{122},$$
$$1 = [a_3^2, a_3] \equiv [c_{13}, a_3] = c_{133},$$
$$1 \equiv c_{23}^2 \equiv [a_3^2, a_2] \equiv c_{132},$$
$$1 \equiv c_{23}^2 c_{232} \equiv [a_2^2, a_3] \equiv c_{123},$$
$$c_{12}^2 \equiv c_{12}^2 c_{122} = [a_1, a_2^2] \equiv c_{121}^{-1} \equiv c_{121},$$
$$c_{13}^2 \equiv c_{13}^2 c_{133} = [a_1, a_3^2] \equiv c_{131} \equiv 1 \pmod{G_4}.$$

Hence,
$$c_{122} \equiv c_{131} \equiv c_{132} \equiv c_{133} \equiv c_{13}^2 \equiv 1, \quad c_{121} \equiv c_{12}^2 \pmod{G_4},$$
whence $G_3 = \langle c_{12}^2, G_4 \rangle$. The arguments as above then show that $G_3 = \langle c_{12}^2 \rangle$. Therefore, $G'$ is abelian with $G' \simeq (2, 2^m)$, as desired. $\square$

**Proposition 12.** *Let $k$ be an imaginary quadratic number field with discriminant $d_k$. Furthermore, let $G = \mathrm{Gal}(k^\infty/k)$. Then $G/G_3 \simeq 32.035$ if and only if there is a factorization of $d_k = -4qq'p$ with $p$, $q$, $q'$ distinct primes satisfying*

(1) $q \equiv q' \equiv 3 \bmod 4$, and $p \equiv 1 \bmod 4$,
(2) $(q/p) = (q'/p) = (q'/q) = -1$, $q \equiv 3 \bmod 8$, $q' \equiv 7 \bmod 8$ and $p \equiv 5 \bmod 8$.

*Moreover, $G \simeq \Gamma_n^{(35)}$, where $n$ is determined by $2^n = h_2(d_3 d_4)$ and $\Gamma_n^{(35)} =$*

$$\langle a_1, a_2, a_3 : \ a_1^4 = a_2^{2^{n+1}} = a_3^4 = 1, \ a_1^2 = c_{13}, \ a_2^2 = c_{12}, \ a_3^2 = c_{13} c_{12}^{2^{n-1}}, \ c_{23} = 1 \rangle.$$

*Proof.* The first part of the proposition follows from the tables. For the structure of the Galois group we refer to [23], Theorem 1. (A group theoretic argument may be given as in the above cases.) $\square$

6.6. **$G/G_3 \simeq 32.036$.**

Using group theory, it can be shown that if $G$ is a finite 2-group such that $G/G_3 \simeq 32.036$, then $G'$ has rank 2. It is, however, not possible to show that $G'' = 1$ because there are finite 2-groups with $G/G_3 \simeq 32.036$ and $G'' \neq 1$ such as the following group $G$ of order $2^9$ for which $G/G_3$ is 32.036 and for which $G''$ is nontrivial:

$$G = \langle a_1, a_2, a_3 : a_1^2 = a_2^2 = c_{12}^{16} = c_{13}^{16} = 1, \ a_3^2 = c_{12}^2 c_{13}^2, \ [c_{ij}, c_{klm}] = 1 \rangle.$$

Next we assume that $k$ is a complex quadratic field with $\mathrm{Cl}_2(k) = (2, 2, 2)$ and $G/G_3 = 32.036$ for $G = \mathrm{Gal}(k^2/k)$; technical problems in an application of the ambiguous class number formula prevents us from proving $G'' = 1$ using only number theory.

The account below combines parts of both approaches. We start with

**Proposition 13.** *Let $G$ be a finite 2-group such that $G/G_3$ is isomorphic to group 32.036. Then $G'$ has rank 2.*

*Proof.* The additional relations are
$$a_1^2 \equiv a_2^2 \equiv 1, \ \ a_3^2 \equiv c_{13} \pmod{G_3}.$$

Thus,
$$1 \equiv [a_1^2, a_j] \equiv c_{1j}^2 c_{1j1},$$
$$1 \equiv [a_2^2, a_j] \equiv c_{2j}^2 c_{2j2},$$
$$1 = [a_3^2, a_3] \equiv [c_{13}, a_3] \equiv c_{133},$$
$$c_{132} \equiv [a_3^2, a_2] \equiv c_{23}^2 \equiv 1 \pmod{G_4},$$

whence,
$$c_{123} \equiv c_{132} \equiv c_{133} \equiv 1, \ \ c_{122} \equiv c_{121} \equiv c_{12}^2, \ \ c_{131} \equiv c_{13}^2 \pmod{G_4}.$$

Arguing as in the proof of Proposition 7, we see that
$$G_j = \langle c_{12}^{2^{j-2}}, c_{13}^{2^{j-2}} \rangle$$
and in particular that $G' = \langle c_{12}, c_{13} \rangle$. □

The following group theoretic lemma will be important below:

**Lemma 3.** *Assume that $G/G_3 = 32.036$ is presented as in the proof of Proposition 13. Let*
$$A = \langle a_2, a_3, G' \rangle \ \text{ and } \ B = \langle a_3, a_1 a_2, G' \rangle.$$
*Then*
$$A'G_4 = \langle c_{23}, c_{12}^2, G_4 \rangle \ \text{ and } \ B'G_4 = \langle c_{13} c_{23}, G_4 \rangle.$$
*Moreover $A/A'$ contains a subgroup of type $(4, 2, 2)$, and $B$ has rank 2.*
*Furthermore, $H = A \cap B$ has rank 2.*
*Finally, the abelianizations of the other five maximal subgroups of $G$ have order 8.*

*Proof.* The first result of the lemma follows immediately by direct calculation. Next notice $A' \subseteq G_3$; but $A/G_3 \simeq (4, 2, 2)$ as can be seen from the presentation of $G/G_3$. Hence $A/A'$ contains a subgroup of type $(4, 2, 2)$. Now let $N$ be the Frattini subgroup of $B$. Then $N = B^2 B' = B^2$ (since $B$ is a 2-group) and $N$ is normal in $G$. Notice
$$\langle c_{13} c_{23}, a_3^2, (a_1 a_2)^2 \rangle \subseteq NG_3 \subseteq G'.$$
But then, since $a_3^2 \equiv c_{13} \pmod{G_3}$ and $(a_1 a_2)^2 \equiv a_1^2 a_2^2 c_{12} \equiv c_{12} \pmod{G_3}$,
$$NG_3 = \langle c_{12}, c_{13} \rangle = G'.$$
Therefore by Theorem 2.49(ii) of [13], $N = G'$. But $B/G'$ has rank 2, whence so does $B$ by the Burnside Basis Theorem, cf. [13] again. Furthermore, $H = \langle a_3, c_{12}, c_{13} \rangle$. But notice that $H^2 = \langle a_3^2, G_3 \rangle$ (recall that $G_3 = \langle c_{12}^2, c_{13}^2 \rangle$). Hence $H/H^2 = \langle \overline{a_3}, \overline{c_{12}} \rangle$, since $a_3^2 \equiv c_{13} \mod G_3$. Since $H/H^2$ has rank 2, so does $H$. Finally, the last statement follows by direct calculation and is left as an exercise to the reader. □

Now assume that $k$ is a complex quadratic field with $\mathrm{Cl}_2(k) = (2,2,2)$ and $G/G_3 = 32.036$ for $G = \mathrm{Gal}(k^2/k)$.

By our table at the end of this article, we see that the discriminant $d_k$ of $k$ is of type $2A, 4C, 4K$, or $4L$. For the rest of this section we reorder the prime discriminants $d_i$ for $i = 1, 2, 3, 4$ according to the following table:

| type | $d_1$ | $d_2$ | $d_3$ | $d_4$ |
|------|-------|-------|-------|-------|
| $2A$ | as | is | | |
| $4C$ | $-p_1$ | $-p_2$ | $-4$ | $p_3$ |
| $4K/L$ | $-4$ | $-p_1$ | $-p_2$ | $p_3$ |

We have the following relations among the Kronecker symbols:
$$(d_1/p_2) = (d_1/p_3) = (d_4/p_1) = (d_2/p_4) = -1,$$
$$(d_2/p_3) = (d_4/p_3) = (d_2/p_1) = 1,$$
where $p_i$ is the prime dividing the new $d_i$.

The following proposition follows immediately from the tables at the end of the article:

**Proposition 14.** *Let $k$ be an imaginary quadratic number field with discriminant $d_k$. Furthermore, let $G = \mathrm{Gal}(k^\infty/k)$. Then $G/G_3 \simeq 32.036$ if and only if there is a factorization of $d_k = d_1 d_2 d_3 d_4$ into distinct prime discriminants satisfying*

(1) $d_i < 0$, $(i = 1, 2, 3)$, and $d_4 > 0$,
(2) $(d_1/p_2) = (d_1/p_3) = (d_4/p_1) = (d_2/p_4) = -1$,
$(d_2/p_3) = (d_4/p_3) = (d_2/p_1) = 1$,
*where $p_i$ is the unique prime dividing $d_i$.*

Now we claim

**Theorem 3.** *Put $M = \mathbb{Q}(\sqrt{d_1 d_2}, \sqrt{d_3}, \sqrt{d_4})$; then $k_{\mathrm{gen}} = M(\sqrt{d_1})$, and we have $h_2(M) = 2^{m+n+1}$, rank $\mathrm{Cl}_2(M) = 2$, and $h_2(k_{\mathrm{gen}}) = 2^{m+n}$. In particular, we have $k^2 = k^3 = M^1$.*

**Remark.** With some more effort, it can be proved that $\mathrm{Cl}_2(M) = (2^{m+1}, 2^n)$ and $\mathrm{Cl}_2(k_{\mathrm{gen}}) = (2^m, 2^n)$, where $2^n = h_2(d_3 d_4)$ and $2^{m+1} = h_2(d_1 d_2 d_3)$.

*Proof.* Since the class numbers of the quadratic subfields of $k_{\mathrm{gen}}$ will occur frequently in our calculations, here's a table for the fields whose class number is even:

| disc $K$ | $d_1 d_4$ | $d_2 d_4$ | $d_3 d_4$ | $d_1 d_2 d_3$ | $d_1 d_2 d_4$ | $d_1 d_3 d_4$ | $d_2 d_3 d_4$ | $d = d_k$ |
|----------|-----------|-----------|-----------|---------------|---------------|---------------|---------------|-----------|
| $h_2(K)$ | $2$ | $2$ | $2^n$ | $2^{m+1}$ | $2$ | $2$ | $2$ | $8$ |

Simple applications of class number formulas show that $k(\sqrt{d_4})$ and $k(\sqrt{d_3 d_4})$ are the only two quadratic unramified extensions of $k$ with 2-class numbers divisible by 16. Thus their compositum $M$ corresponds to the Galois group $A \cap B$ in Lemma 3, and we conclude that $\mathrm{Cl}_2(M)$ has rank 2.

Now let us compute the class number of $M$. We first determine the unit index $q(M/\mathbb{Q})$, and to do so we have to show that $q(M^+/\mathbb{Q}) = 2$, where $M^+ = \mathbb{Q}(\sqrt{d_1 d_2}, \sqrt{d_4})$ is the maximal real subfield of $M$. Since $M^+$ is an unramified quadratic extension of the field $F = \mathbb{Q}(\sqrt{d_1 d_2 d_4})$ with 2-class number 2, $M^+$ is the 2-class field of $F$ and has odd class number (since $\text{Cl}_2(F)$ is cyclic). On the other hand, the class number formula (1) says

$$\begin{aligned} h_2(M^+) &= 2^{-2}q(M^+/\mathbb{Q})h_2(d_1d_2)h_2(d_4)h_2(d_1d_2d_4) \\ &= 2^{-2}q(M^+/\mathbb{Q}) \cdot 2. \end{aligned}$$

Thus $q(M^+/\mathbb{Q}) = 2$.

Now $q(M/\mathbb{Q}) = Q(M)q(M^+/\mathbb{Q})$, hence we can apply the class number formula to $M/\mathbb{Q}$ once we know $Q(M)$. If $d_3$ is odd, then $M/M^+$ is essentially ramified, hence $Q(M) = 1$ by [21], p. 350, Theorem 1. If $d_3 = -4$, then we are in case (ii).1. of [21], Theorem 1, so again $Q(M) = 1$, hence $q(M/\mathbb{Q}) = Q(M)q(M^+/\mathbb{Q}) = 2$. Now

$$h_2(M) = 2^{-5}q(M/\mathbb{Q})h_2(d_3d_4)h_2(d_1d_2d_3)h_2(d_1d_2d_4)h_2(k) = 2^{m+n+1}.$$

Now let us compute $h_2(K)$, where $K = k_{\text{gen}}$. We first claim that $K^+$ has odd class number:

**Lemma 4.** *Let $k$ be an imaginary quadratic field and let $G = \text{Gal}(k^2/k)$. Suppose furthermore that $G/G_3$ is isomorphic to group 32.036. Then the maximal real subfield of $k^1$ has odd class number.*

*Proof.* We assume $d_k = d_1d_2d_3d_4$, where the $d_j$ are ordered as in the table given above. Let $p_j$ be the primes dividing these new $d_j$'s. Let $K = k^1 = k_{\text{gen}} = \mathbb{Q}(\sqrt{d_1}, \sqrt{d_2}, \sqrt{d_3}, \sqrt{d_4})$ and $K^+ = K \cap \mathbb{R} = \mathbb{Q}(\sqrt{d_1d_2}, \sqrt{d_1d_3}, \sqrt{d_4})$. We need to show that $h_2(K^+) = 1$. To this end, consider the subfield $L = \mathbb{Q}(\sqrt{d_1d_3}, \sqrt{d_4})$ of $K^+$. Since $L$ is the Hilbert 2-class field of the field $k = \mathbb{Q}(\sqrt{d_1d_3d_4})$ with 2-class number 2, we deduce that $h(L)$ is odd. Now we apply the ambiguous class number formula to $K^+/L$: since exactly the two prime ideals above $p_2$ ramify, we have $t = 2$, hence $\#\text{Am}(K^+/L) = 2/(E:H)$, and our claim will follow if we can exhibit a unit in $E = E_L$ that is not a norm from $K^+$. By Hasse's norm theorem, a unit is a norm from a quadratic extension if and only if it is a norm in all local extensions $K^+_\mathfrak{P}/L_\mathfrak{p}$; since units are always norms in unramified extensions of local fields, we only have to study the localizations at the primes above $p_2$. The lemma below shows that the fundamental unit $\varepsilon_4$ of $\mathbb{Q}(\sqrt{d_4})$ is not a quadratic residue modulo the primes above $p_2$, hence it is not a local norm at these primes by Hensel's lemma, and we are done.  $\square$

**Lemma 5.** *Let $K$ be a number field containing a real quadratic field $k$ whose fundamental unit $\varepsilon$ has negative norm. If $p \equiv 3 \bmod 4$ is a prime that is inert in $k$ and splits completely in $K/k$, then $\varepsilon$ is a quadratic nonresidue modulo the primes above $p$.*

*Proof.* Let $\mathfrak{p}$ be a prime ideal above $p$ in $\mathcal{O}_K$; by assumption it has absolute norm $p^2$, hence $[\varepsilon/\mathfrak{p}] \equiv \varepsilon^{(p^2-1)/2} = (\varepsilon^{p+1})^{(p-1)/2} \bmod \mathfrak{p}$. Now $\varepsilon^{p+1} \equiv N_{k/\mathbb{Q}}\varepsilon = -1 \bmod p$ for elementary reasons, hence $[\varepsilon/\mathfrak{p}] \equiv (-1)^{(p-1)/2} = -1 \bmod \mathfrak{p}$.  $\square$

The fact that $2 \nmid h(K^+)$ allows us to show that $q(K^+/\mathbb{Q}) = 2^6$. In fact, the class number formula gives

$$h_2(K^+) = 2^{-9}q(K^+/\mathbb{Q})\prod h_r$$

with $\prod h_r = h_2(d_4)h_2(d_1d_2)h_2(d_1d_3)h_2(d_2d_3)h_2(d_1d_2d_4)h_2(d_1d_3d_4)h_2(d_2d_3d_4)$. If we plug in the class numbers from the table above we get $h_2(K^+) = 2^{-6}q(K^+/\mathbb{Q})$. Since $h_2(K^+) = 1$, the claim follows.

Now we apply the class number formula to $K/\mathbb{Q}$. We have $Q(K) = 2$ by [21], p. 352, Example 4, hence $h_2(K) = 2^{-16} \cdot 2^7 \prod h_i$, where the table above shows $\prod h_i = 2^8 h_2(d_3 d_4) h_2(d_1 d_2 d_3)$, i.e., $h_2(K) = 2^{m+n}$. By [3], Proposition 7, we conclude that $k^2 = k^3$. □

Let us consider the cases 4K/L, that is we assume that $d = -4pqq'$ with $p \equiv 5 \bmod 8$, $q \equiv 7 \bmod 8$, $(p/q) = -1$, $(p/q') = +1$, $(q/q') = -1$. Put $2^n = h_2(-pq')$ and $2^m = h_2(-qq')$.

| $j$ | $k_j$ | $\kappa_j$ | $h(k_j)$ | rank $C_j$ | $\mathrm{Cl}_2(k_j)$ | $NC_j$ | |
|---|---|---|---|---|---|---|---|
| 1 | $k(\sqrt{p})$ | $\langle[\mathfrak{q}\mathfrak{q}']\rangle$ | $2^{m+2}$ | 3 | $(2,2,2^m)$ | $\langle[2\mathfrak{q}],[\mathfrak{q}']\rangle$ | |
| 2 | $k(\sqrt{-q})$ | $\langle[\mathfrak{q}],[\mathfrak{q}']\rangle$ | 8 | 3 | $(2,2,2)$ | $\langle[2],[\mathfrak{q}']\rangle$ | |
| 3 | $k(\sqrt{q'})$ | $\langle[2],[\mathfrak{q}']\rangle$ | 8 | 3 | $(2,2,2)$ | $\langle[\mathfrak{q}],[\mathfrak{q}']\rangle$ | |
| 4 | $k(\sqrt{-1})$ | $\langle[2],[\mathfrak{q}']\rangle$ | 8 | 2 | $(2,4)$ | $\langle[2],[\mathfrak{q}\mathfrak{q}']\rangle$ | $\langle[2\mathfrak{q}],[\mathfrak{q}\mathfrak{q}']\rangle$ |
| 5 | $k(\sqrt{-q'})$ | $\langle[2\mathfrak{q}],[\mathfrak{q}']\rangle$ | 8 | 3 | $(2,2,2)$ | $\langle[2\mathfrak{q}],[\mathfrak{q}\mathfrak{q}']\rangle$ | $\langle[2],[\mathfrak{q}\mathfrak{q}']\rangle$ |
| 6 | $k(\sqrt{-p})$ | $\langle[\mathfrak{q}],[\mathfrak{q}']\rangle$ | 8 | 2 | $(2,4)$ | $\langle[2\mathfrak{q}'],[\mathfrak{q}]\rangle$ | $\langle[2],[\mathfrak{q}]\rangle$ |
| 7 | $k(\sqrt{q})$ | $\langle[2],[\mathfrak{q}]\rangle$ | $2^{n+2}$ | 2 | $(2,2^{n+1})$ | $\langle[2],[\mathfrak{q}]\rangle$ | $\langle[2\mathfrak{q}'],[\mathfrak{q}]\rangle$ |

Again, the left column in $NC_j$ is for primes $q' \equiv 3 \bmod 8$, the right one for $q' \equiv 7 \bmod 8$.

Numerical examples:

| $d$ | $p$ | $q$ | $q'$ | $m$ | $n$ |
|---|---|---|---|---|---|
| $-1540$ | 5 | 7 | 11 | 3 | 2 |
| $-7332$ | 13 | 47 | 3 | 3 | 2 |
| $-8372$ | 13 | 7 | 23 | 4 | 3 |
| $-10212$ | 37 | 23 | 3 | 3 | 2 |

### 6.7. $G/G_3 \simeq 32.033$.

If $G = \mathrm{Gal}(k^2/k)$ with $k$ imaginary quadratic and $G/G_3$ isomorphic to group 32.033, then we shall see that $\mathrm{rank}\,\mathrm{Cl}_2(k^1) \geq 3$. It is not hard to find many examples of $k$ with infinite 2-class field tower, for example, anytime the 2-class number of $\mathbb{Q}(\sqrt{d_1 d_2 d_4})$ for $d_k$ of type 2B and $\mathbb{Q}(\sqrt{-d_1 d_3})$ for the other types is $\geq 16$; compare with [27], Beispiel 4.

**Proposition 15.** *Suppose $k$ is an imaginary quadratic field and $G = \mathrm{Gal}(k^2/k)$ is such that $G/G_3 \simeq 32.033$. Then $\mathrm{rank}(G') \geq 3$.*

*Proof.* This is an immediate consequence of the fact that the Schur multiplier of group 32.033 is $(2,2,2,2)$ and so in particular of rank 4. If $K$ denotes the fixed field of $G_3$ in $k^2$, then

$$\begin{aligned}
\mathrm{rank}\,\mathrm{Gal}(k^2/k^1) &= \mathrm{rank}(G'/G'') \geq \mathrm{rank}(G_3/G_4) \\
&= \mathrm{rank}\,\mathrm{Gal}(K_{\mathrm{cen}}/K_{\mathrm{gen}}) \\
&\geq \mathrm{rank}\,\mathcal{M}(G/G_3) - \mathrm{rank}\,E_k/(E_k \cap N_{K/k}(K^\times)) \geq 3,
\end{aligned}$$

where $K_{\text{cen}}$ and $K_{\text{gen}}$ are the central and genus class field extensions of $K/k$, respectively, (which are in our case the fixed fields of $G_4$ and $G_3$ ($K_{\text{gen}} = K$)) and $\mathcal{M}$ is the Schur multiplier, see e.g. [2] for more information. $\square$

From the tables, we have

**Proposition 16.** *Let $k$ be an imaginary quadratic number field with discriminant $d_k$. Furthermore, let $G = \text{Gal}(k^\infty/k)$. Then $G/G_3 \simeq 32.033$ if and only if there is a factorization of $d_k = d_1 d_2 d_3 d_4$ into distinct prime discriminants satisfying*

(1) $d_i < 0$, $(i = 1, 2, 3)$, and $d_4 > 0$,
(2) $(d_1/p_2) = (d_2/p_3) = (d_3/p_1) = (d_1/p_4) = -1$, $(d_4/p_2) = (d_4/p_3) = 1$, *where $p_i$ is the unique prime dividing $d_i$.*

Now consider the unramified quadratic extension $K = k(\sqrt{d_4})$ of $k$, where $k$ satisfies the conditions in Proposition 16. The class number formula immediately gives $h_2(K) = 16$, and the ambiguous class number formula applied to $K/\mathbb{Q}(\sqrt{p})$ tells us that $\text{Cl}_2(K)$ has 2-rank 4. Thus $\text{Cl}_2(K) \simeq (2,2,2,2)$.

Now we make two claims:

(1) $\text{Cl}_2(k^1)$ has rank 3;
(2) $\text{Cl}_2(K^1)$ has rank at least 4.

It is an immediate consequence of these assertions that we must have $k^2 \neq k^3$: if we had $k^2 = k^3$, then the ranks of the class groups of the intermediate fields of $k^2/k^1$ could not increase.

Ad 1.: We have already shown that $\text{Cl}_2(k^1)$ has rank at least 3; to show that it is at most 3, we apply the ambiguous class number formula to $k^1/L$, where $L = \mathbb{Q}(\sqrt{d_1}, \sqrt{d_2}, \sqrt{d_3})$. Since $L$ is a quadratic unramified extension of $F = \mathbb{Q}(\sqrt{d_1}, \sqrt{d_2 d_3})$, and since moreover $h_2(F) = 2$ by the class number formula, we deduce that $L$ has odd class number. Since there are exactly 4 primes above $p_4$ in $L$, the ambiguous class number formula says that $\text{rank Cl}_2(L) \leq 3$ as desired.

Ad 2.: Let $M = K^1$ with $K = k(\sqrt{d_4})$. Then $H = \text{Gal}(M/K) \simeq (2,2,2,2)$ has Schur multiplier of rank $\binom{4}{2} = 6$, thus $\text{rank Gal}(M_{\text{cen}}/M_{\text{gen}}) \geq 6 - \text{rank } E_K/E_K^2 \geq 4$.

This also concludes our proof of Theorem 1.

Now consider the special case $d = -4pqq'$, where $p \equiv 1 \bmod 8$, $q \equiv 3 \bmod 4$ and $q' \equiv 3 \bmod 8$ are primes such that $(p/q) = -1$, $(p/q') = +1$ and $(q/q') = +1$; we are in case 4E if $q \equiv 3 \bmod 8$ and in case 4F if $q \equiv 7 \bmod 8$. Define integers $l, m, n \geq 2$ by $2^l = h_2(-pq')$, $2^m = h_2(4pq')$ and $2^n = h_2(-4p)$, where $h_2(D)$ denotes the 2-class number of the quadratic field with discriminant $D$.

| $j$ | $k_j$ | $h(k_j)$ | rank $C_j$ | $\text{Cl}_2(k_j)$ | $NC_j$ | |
|---|---|---|---|---|---|---|
| 1 | $k(\sqrt{p})$ | 16 | 4 | $(2,2,2,2)$ | $\langle [2], [\mathfrak{p}\mathfrak{q}] \rangle$ | |
| 2 | $k(\sqrt{-q'})$ | 8 | 3 | $(2,2,2)$ | $\langle [\mathfrak{p}], [\mathfrak{q}] \rangle$ | |
| 3 | $k(\sqrt{q})$ | $2^{l+2}$ | 2 | $(2, 2^{l+1})$ | $\langle [\mathfrak{p}\mathfrak{q}], [2\mathfrak{p}] \rangle$ | |
| 4 | $k(\sqrt{-1})$ | 8 | 3 | $(2,2,2)$ | $\langle [\mathfrak{p}], [2] \rangle$ | $\langle [\mathfrak{p}], [2\mathfrak{q}] \rangle$ |
| 5 | $k(\sqrt{-p})$ | $2^{n+2}$ | 2 | $(2, 2^{n+1})$ | $\langle [\mathfrak{q}], [2] \rangle$ | $\langle [\mathfrak{q}], [2\mathfrak{p}] \rangle$ |
| 6 | $k(\sqrt{-q})$ | $2^{m+2}$ | 2 | $(2, 2^{m+1})$ | $\langle [\mathfrak{q}], [2\mathfrak{p}] \rangle$ | $\langle [\mathfrak{q}], [2] \rangle$ |
| 7 | $k(\sqrt{q'})$ | 8 | 3 | $(2,2,2)$ | $\langle [\mathfrak{p}], [2\mathfrak{q}] \rangle$ | $\langle [\mathfrak{p}], [2] \rangle$ |

Since $\mathrm{Gal}(k_2^1/k)$ is a group of order 32 with abelianization $(2,2,2)$ and a subgroup of type $(2,2,2,2)$, we deduce that $\mathrm{Gal}(k_2^1/k) \simeq 32.033$.

Let us show that $\mathrm{Cl}_2(k_5)$ has rank 2. We apply the ambiguous class number formula to $k_5/\mathbb{Q}(\sqrt{qq'})$; there are exactly four ramified places, and $-1$ is not a local norm at the infinite primes. Thus our claim is equivalent to showing that $\varepsilon_{qq'}$ is a local norm everywhere. We first observe that $q\varepsilon_{qq'}$ is a square, and $(-p/q) = +1$ guarantees that $q$ is a local norm everywhere except possibly at the prime 2: but the product formula takes care of that.

Numerical examples:

|     $d$ | case |   $p$ | $q$ | $q'$ | $l$ | $m$ | $n$ |
|--------:|:----:|------:|:---:|:----:|:---:|:---:|:---:|
|  $-6132$ |  4F  |   73  |  7  |   3  |  2  |  2  |  2  |
|  $-8148$ |  4F  |   97  |  7  |   3  |  2  |  2  |  2  |
| $-11748$ |  4E  |   89  |  3  |  11  |  3  |  2  |  2  |
| $-12036$ |  4E  |   17  |  3  |  59  |  2  |  2  |  2  |
| $-14916$ |  4E  |  113  |  3  |  11  |  2  |  2  |  3  |
| $-26292$ |  4F  |  313  |  7  |   3  |  3  |  2  |  3  |
| $-40836$ |  4E  |   41  |  3  |  83  |  3  |  3  |  3  |

## Appendix A

In this appendix we indicate how one goes about constructing some of the unramified extensions explicitly. Let us recall the following lemma taken from [22]:

**Lemma 6.** *Let $K/F$ be a quartic extension with $\mathrm{Gal}(K/F) \simeq (2,2)$; let $\sigma, \tau$ and $\sigma\tau$ denote its nontrivial automorphisms, and put $M = K(\sqrt{\mu})$. Then $M/F$ is normal if and only if $\mu^{1-\rho} \stackrel{2}{=} 1$ for all $\rho \in \mathrm{Gal}(K/F)$. If this is the case, write $\mu^{1-\sigma} = \alpha_\sigma^2$, $\mu^{1-\tau} = \alpha_\tau^2$ and $\mu^{1-\sigma\tau} = \alpha_{\sigma\tau}^2$. It is easy to see that $\alpha_\rho^{1+\rho} = \pm 1$ for all $\rho \in \mathrm{Gal}(K/F)$; define $S(\mu, K/F) = (\alpha_\sigma^{1+\sigma}, \alpha_\tau^{1+\tau}, \alpha_{\sigma\tau}^{1+\sigma\tau})$ and identify vectors which coincide upon permutation of their entries. Then*

$$\mathrm{Gal}(M/F) \simeq \begin{cases} (2,2,2) & \iff S(\mu, K/F) = (+1,+1,+1), \\ (2,4) & \iff S(\mu, K/F) = (-1,-1,+1), \\ D_4 & \iff S(\mu, K/F) = (-1,+1,+1), \\ H_8 & \iff S(\mu, K/F) = (-1,-1,-1). \end{cases}$$

*Moreover, $M$ is cyclic over the fixed field of $\langle \rho \rangle$ if and only if $\alpha_\rho^{1+\rho} = -1$, and has type $(2,2)$ otherwise.*

In the following, we show how to find the 2-class fields for $G = 32.041$ if $d = \mathrm{disc}\, k$ is odd.

Let the prime discriminants $d_i$ be as in Proposition 4. Then the following equations have integral solutions by Legendre's theorem:

$$\begin{align}
(2) \qquad d_1 x^2 - d_4 y^2 &= d_2 z^2, \\
(3) \qquad d_2 x^2 - d_4 y^2 &= d_3 z^2, \\
(4) \qquad d_3 x^2 - d_4 y^2 &= d_1 z^2.
\end{align}$$

The same is true for the following equations:

$$
\begin{align}
t^2 - d_2 d_3 u^2 &= -d_4 v^2, \tag{5}\\
t^2 - d_3 d_1 u^2 &= -d_4 v^2, \tag{6}\\
t^2 - d_1 d_2 u^2 &= -d_4 v^2. \tag{7}
\end{align}
$$

Let $(x_1, y_1, z_1) \in \mathbb{Z} \times \mathbb{Z} \times \mathbb{Z}$ be a solution of (2) such that $(x_1, y_1) = (x_1, z_1) = (y_1, z_1) = 1$; we will call such solutions primitive. It follows from Cassels [10] that there are primitive solutions such that $2 \mid y_1$ and $2 \nmid x_1 z_1$. Similarly, we may demand that, say, equation (5) has a solution $(t_1, u_1, v_1)$ with $t_1$ even and $u_1$ odd. Now put

$$
\begin{align}
\alpha_1 &= x_1 \sqrt{d_1} + y_1 \sqrt{d_4}, \\
\alpha_2 &= x_2 \sqrt{d_2} + y_2 \sqrt{d_4}, \\
\beta_3 &= t_3 + u_3 \sqrt{d_1 d_2}
\end{align}
$$

and form the product $\mu_1 = \alpha_1 \alpha_2 \beta_3$. Then $\mu_1 \equiv (\sqrt{d_1}+y_1)(\sqrt{d_2}+y_2)(t_3+\sqrt{d_1 d_2}) \equiv (y_1 + y_2 + \sqrt{d_1 d_2})(t_3 + \sqrt{d_1 d_2}) \equiv 1 + y_1 + y_2 + t_3 \bmod 4$. Replacing $\mu_1$ by $-\mu_1$ if $y_1 + y_2 + t_3 \equiv 2 \bmod 4$ we find that $\mu_1 \equiv 1 \bmod 4$. Since $(\mu_1)$ is an ideal square, $K(\sqrt{\mu})/K$ is an unramified extension of $K = \mathbb{Q}(\sqrt{d_1}, \sqrt{d_2}, \sqrt{d_3}, \sqrt{d_4})$. The first thing to do is check that this extension is quadratic, i.e. that $\mu$ is not a square in $K$. This is not obvious, as the following observation shows:

**Lemma 7.** *Put $\gamma := \pm \beta_1 \beta_2 \beta_3$ and choose the sign so that $\gamma$ is a square modulo $4$. Then $\gamma \stackrel{2}{=} d_4$ in $K_0 = \mathbb{Q}(\sqrt{d_1}, \sqrt{d_2}, \sqrt{d_3})$.*

*Proof.* (Sketch) Put $F = \mathbb{Q}(\sqrt{d_1 d_2}, \sqrt{d_1 d_3}, \sqrt{d_2 d_3})$. Using Lemma 6, we easily show that $F(\sqrt{\gamma})/\mathbb{Q}$ is elementary abelian. Since $F(\sqrt{\gamma})/F$ is unramified outside $d_4 \infty$, this implies that $\gamma \stackrel{2}{=} d_4 m$ for some discriminant $m \mid d_1 d_2 d_3$. The claim follows, since $m$ becomes a square in $K_0$. $\square$

Thus our candidates for unramified quadratic extensions of $K$ are $\mu = \alpha_1 \alpha_2 \beta_3$ and $\nu = \alpha_2 \alpha_3 \beta_1$. Note that $\mu\nu \stackrel{2}{=} \alpha_1 \alpha_3 \beta_1 \beta_3 \stackrel{2}{=} \alpha_1 \alpha_3 \beta_2$ in $K$. Our next job is to compute the Galois group $G_M = \mathrm{Gal}(M/k)$, where $M = K(\sqrt{\mu})$. Note that, by symmetry, the Galois groups of $K(\sqrt{\nu})/k$ and $K(\sqrt{\mu\nu})/k$ must be isomorphic to $G_M$. Since these are subgroups of index 2 in the group $\Gamma = \mathrm{Gal}(k^2/k)$ of order 32, and since we know that $\Gamma' \simeq (2,2)$ and $\Gamma/\Gamma' \simeq (2,2,2)$, the group tables in [14] already imply that $\Gamma \simeq 32.041$ as soon as we know that $M \neq K$ (note that this shows that $M \neq N = K(\sqrt{\nu})$, since $M = N$ implies that $\mu\nu$ is a square in $K$ which it isn't). This will be done by computing some Galois groups.

In fact, we put $K_i = k(\sqrt{d_i})$ and then show that $\mathrm{Gal}(M/K_1) \simeq H_8$ and $\mathrm{Gal}(M/K_2) \simeq D_4$. Thus $G_M$ is a group of order 16 with $G'_M \simeq \mathbb{Z}/2\mathbb{Z}$, $G_M/G'_M \simeq (2,2,2)$, and possessing subgroups of type $H_8$ and $D_4$. The only such group is $D_4 \curlyvee C_4$. Here are the details:

- $\mathrm{Gal}(M/K_1) \simeq H_8$. Define automorphisms $\sigma, \tau \in \mathrm{Gal}(K/K_1)$ by

|   | $\sqrt{d_1}$ | $\sqrt{d_2}$ | $\sqrt{d_3}$ | $\sqrt{d_4}$ |
|---|---|---|---|---|
| $\sigma$ | $+\sqrt{d_1}$ | $-\sqrt{d_2}$ | $-\sqrt{d_3}$ | $+\sqrt{d_4}$ |
| $\tau$ | $+\sqrt{d_1}$ | $-\sqrt{d_2}$ | $+\sqrt{d_3}$ | $-\sqrt{d_4}$ |
| $\sigma\tau$ | $+\sqrt{d_1}$ | $+\sqrt{d_2}$ | $-\sqrt{d_3}$ | $-\sqrt{d_4}$ |

Then we find $\alpha_\sigma = \alpha_2\beta_3/\sqrt{d_3d_4}z_2^2v_3^2$, $\alpha_\tau = \alpha_1\beta_3/\sqrt{d_2d_4}z_1^2v_3^2$, as well as $\alpha_{\sigma\tau} = \alpha_1\alpha_2/\sqrt{d_2d_3}z_1^2z_2^2$, hence $\alpha_\sigma^{1+\sigma} = \alpha_\tau^{1+\tau} = \alpha_{\sigma\tau}^{1+\sigma\tau} = -1$, and our claim follows.

- $\mathrm{Gal}(M/K_2) \simeq D_4$. As above we start by defining some automorphisms of $\mathrm{Gal}(K/K_2)$:

|   | $\sqrt{d_1}$ | $\sqrt{d_2}$ | $\sqrt{d_3}$ | $\sqrt{d_4}$ |
|---|---|---|---|---|
| $\sigma$ | $-\sqrt{d_1}$ | $+\sqrt{d_2}$ | $-\sqrt{d_3}$ | $+\sqrt{d_4}$ |
| $\tau$ | $-\sqrt{d_1}$ | $+\sqrt{d_2}$ | $+\sqrt{d_3}$ | $-\sqrt{d_4}$ |
| $\sigma\tau$ | $+\sqrt{d_1}$ | $+\sqrt{d_2}$ | $-\sqrt{d_3}$ | $-\sqrt{d_4}$ |

Since $\sigma\tau$ is the "same" automorphism as in the case discussed above, we get $\alpha_{\sigma\tau}^{1+\sigma\tau} = -1$. Next, $\alpha_\sigma = \alpha_1\beta_3/\sqrt{d_2d_4}z_1^2v_3^2$, hence $\alpha_\sigma^{1+\sigma} = +1$. Similarly, $\alpha_\tau = \alpha_2\beta_3/\sqrt{d_3d_4}z_2^2v_3^2$ and $\alpha_\tau^{1+\tau} = +1$.

Example: Take $d_1 = -7$, $d_2 = -3$, $d_3 = -23$ and $d_4 = 5$. Then $\mathrm{Cl}(k) \simeq (2, 2, 2, 5)$. Here $\alpha_1 = 2\sqrt{5} + \sqrt{-7}$, $\alpha_2 = 2\sqrt{5} + \sqrt{-3}$, and $\beta_3 = 4 + \sqrt{21}$; then $\mu = -\alpha_1\alpha_2\beta_3$ is primary, and the field $L_\mu = \mathbb{Q}(\sqrt{\mu})$ is also generated by a root of the polynomial $h(x) = x^8 - 3x^7 + 14x^6 - 38x^5 + 85x^4 - 134x^3 + 244x^2 - 120x + 240$. The field discriminant is disc $L_\mu = 3^4 5^4 7^4 23^2$, and pari computes its class group as $\mathrm{Cl}(L_\mu) = \mathbb{Z}/2\mathbb{Z}$.

## Appendix B.

We take this opportunity to give a simple proof of a generalization of a result due to the first author and Parry [5]:

**Theorem 4.** *Let $k$ be a complex quadratic number field with discriminant $d$, and assume that $d = d_1d_2d_3d_4$ is the product of three positive ($d_1$, $d_2$, $d_3$) and one negative ($d_4$) prime discriminants. Then the $2$-class number of the Hilbert $2$-class field $k^1$ of $k$ is divisible by $2^5$.*

*Proof.* Let $L = \mathbb{Q}(\sqrt{d_1}, \sqrt{d_2}, \sqrt{d_3}, \sqrt{d_4})$ denote the genus class field of $k$, and let $K = \mathbb{Q}(\sqrt{d_1}, \sqrt{d_2}, \sqrt{d_3})$ be its maximal real subfield. Since $Q(L) = 1$ by [21], we have $q(L/\mathbb{Q}) = q(K/\mathbb{Q})$. The class number formula (1) applied to $L/\mathbb{Q}$ and $K/\mathbb{Q}$ gives

$$(8) \qquad h_2(L) = 2^{-16}q(L/\mathbb{Q})\prod_j h_j, \quad h_2(K) = 2^{-9}q(K/\mathbb{Q})\prod_r h_r,$$

where the products are over all quadratic subfields $k_j$ and $k_r$ of $L/\mathbb{Q}$ and $K/\mathbb{Q}$, respectively, and where $h_j$ and $h_r$ denote the 2-class numbers of $k_j$ and $k_r$. Thus $h_2(L)/h_2(K) = 2^{-7}\prod_c h_c$, where the product is over all complex quadratic subfields. Let $2^m$ be the 2-class number of $k$; the three fields whose discriminants are products of three (resp. two) prime discriminants have class number divisible by 4 (resp. 2), hence $2^{m+9} \mid \prod_c h_c$. This implies that $2^{m+2} \mid h_2(L)$.

Since $k^1/L$ is an unramified abelian extension of relative degree $2^{m-3}$, we have $2^{3-m} h_2(L) \mid h_2(k^1)$, and this finally proves that $h_2(k^1) \equiv 0 \bmod 2^5$. $\square$

$d_k$ for which $\mathrm{Cl}_2(k) \simeq (2,2,2)$

| Case | Graph | Presentation of $G/G_3$ | Group Number |
|------|-------|--------------------------|--------------|
| 1A | | $s_1^2 = s_2^2 = s_3^2 = 1$ | 64.144 |
| 1B | | $s_1^2 = s_2^2 = 1,\ s_3^2 = t_{23}$ | 64.144 |
| 1C | | $s_1^2 = 1,\ s_2^2 = t_{12},\ s_3^2 = t_{13}$ | 64.144 |
| 1D | | $s_1^2 = 1,\ s_2^2 = t_{12},\ s_3^2 = t_{12}t_{23}$ | 64.144 |
| 1E | | $s_1^2 = s_2^2 = t_{12},\ s_3^2 = 1$ | 64.147 |
| 1F | | $s_1^2 = s_2^2 = t_{12},\ s_3^2 = t_{12}t_{13},$ | 64.147 |
| 1G | | $s_1^2 = 1,\ s_2^2 = t_{13}t_{23},\ s_3^2 = t_{12}t_{23}$ | 64.147 |
| 2A | | $s_1^2 = s_2^2 = 1,\ s_3^2 = t_{13},$ $t_{12}t_{13}t_{23} = 1$ | 32.036 |
| 2B | | $s_1^2 = 1,\ s_2^2 = t_{12},\ s_3^2 = t_{13},$ $t_{12}t_{13}t_{23} = 1$ | 32.033 |
| 2C | | $s_1^2 = s_2^2 = t_{12},\ s_3^2 = t_{13},$ $t_{12}t_{13}t_{23} = 1$ | 32.037 |
| 2D | | $s_1^2 = t_{12},\ s_2^2 = t_{23},\ s_3^2 = t_{13},$ $t_{12}t_{13}t_{23} = 1$ | 32.041 |

$d_k$ for which $Cl_2(k) \simeq (2,2,2)$ (cont.)

| Case | Graph | Presentation of $G/G_3$ | Group Number |
|---|---|---|---|
| 3A | | $s_1^2 = t_{12}t_{23},\ s_2^2 = t_{12}t_{13},\ s_3^2 = 1$ | 64.147 |
| 3B | | $s_1^2 = t_{12}t_{13},\ s_2^2 = t_{12}t_{23},\ s_3^2 = t_{13}t_{23}$ | 64.145 |
| 3C | | $s_1^2 = t_{12}t_{13},\ s_2^2 = t_{13},\ s_3^2 = t_{12}$ | 64.146 |
| 3D | | $s_1^2 = t_{13},\ s_2^2 = t_{12}t_{13},\ s_3^2 = 1$ | 64.144 |
| 3E | | $s_1^2 = 1,\ s_2^2 = t_{12},\ s_3^2 = t_{13}$ | 64.144 |
| 4A | | $s_1^2 = 1,\ s_2^2 = t_{13}t_{23},\ s_3^2 = t_{23}$, $t_{12} = 1$ | 32.037 |
| 4B | | $s_1^2 = t_{23},\ s_2^2 = t_{13},\ s_3^2 = t_{23}$, $t_{12} = 1$ | 32.040 |
| 4C | | $s_1^2 = 1,\ s_2^2 = t_{13}t_{23},\ s_3^2 = t_{13}$, $t_{12} = 1$ | 32.036 |
| 4D | | $s_1^2 = t_{13},\ s_2^2 = t_{23},\ s_3^2 = t_{13}$, $t_{12} = 1$ | 32.035 |
| 4E | | $s_1^2 = s_2^2 = 1,\ s_3^2 = t_{23}$, $t_{12} = 1$ | 32.033 |
| 4F | | $s_1^2 = s_2^2 = 1,\ s_3^2 = t_{23}$, $t_{12} = 1$ | 32.033 |

$d_k$ for which $Cl_2(k) \simeq (2,2,2)$ (cont.)

| Case | Graph | Presentation of $G/G_3$ | Group Number |
|------|-------|--------------------------|--------------|
| 4G | | $s_1^2 = t_{23},\ s_2^2 = 1,\ s_3^2 = t_{23},$ $t_{12} = 1$ | 32.038 |
| 4H | | $s_1^2 = t_{23},\ s_2^2 = 1,\ s_3^2 = t_{23},$ $t_{12} = 1$ | 32.038 |
| 4I | | $s_1^2 = s_2^2 = 1,\ s_3^2 = t_{13},$ $t_{12} = 1$ | 32.033 |
| 4J | | $s_1^2 = s_2^2 = 1,\ s_3^2 = t_{13},$ $t_{12} = 1$ | 32.033 |
| 4K | | $s_1^2 = t_{13},\ s_2^2 = s_3^2 = 1,$ $t_{12} = 1$ | 32.036 |
| 4L | | $s_1^2 = t_{13},\ s_2^2 = s_3^2 = 1,$ $t_{12} = 1$ | 32.036 |

Addresses of the authors:
E.B.
Department of Mathematics
Unity College, Unity, ME 04988 USA
*e-mail:* benjamin@mint.net

F.L.
CSU San Marcos, Dept. Mathematics
San Marcos, CA 92096-0001
*e-mail:* franzl@csusm.edu

C.S.
Department of Mathematics and Statistics
University of Maine
Orono, ME 04469 USA
and
Research Institute of Mathematics, Orono
*e-mail:* snyder@gauss.umemat.maine.edu